\documentclass[12pt,british]{amsart}

\usepackage[latin1]{inputenc}
\usepackage[T1]{fontenc}
\usepackage{babel,eucal,url,amssymb,stmaryrd,booktabs,%
enumerate,amscd,paralist,amsmath,yhmath}

\theoremstyle{plain}
\newtheorem{theorem}{Theorem}[section]

\newtheorem{lemma}[theorem]{Lemma}

\theoremstyle{definition}
\newtheorem{definition}[theorem]{Definition}

\theoremstyle{remark}

\numberwithin{equation}{section}

\newcommand{\norm}[1]{\left\lVert #1\right\rVert}

\newcommand{\cd}{\mathcal B}

\newenvironment{spmatrix}{\left(\smallmatrix}{\endsmallmatrix\right)}

\begin{document}

\title[The quasi-nilpotent DT-operator]{Moment formulas for \\ 
the quasi-nilpotent DT-operator} 

\author{Lars Aagaard}
\address[Lars Aagaard]{Department of Mathematics and Computer Science\\
University of Southern Denmark\\
Campusvej 55\\
DK-5230 Odense M\\
Denmark}
\email{laa@imada.sdu.dk}

\author{Uffe Haagerup$^\dagger$}
\address[Uffe Haagerup]{Department of Mathematics and Computer Science\\
University of Southern Denmark\\
Campusvej 55\\
DK-5230 Odense M\\
Denmark}
\email{Haagerup@imada.sdu.dk}

\begin{abstract}
Let $T$ be the quasi-nilpotent DT-operator. By use of Voiculescu's amalgamated $R$-transform we compute the momets of $(T-\lambda 1)^*(T-\lambda 1)$ where $\lambda\in \mathbb C$, and the Brown-measure of $T+\sqrt{\epsilon}Y$, where $Y$ is a circular element $*$-free from $T$ for $\epsilon>0$. Moreover we give a new proof of \'Sniady's formula for the moments $\tau(((T^*)^k T^k)^n)$ for $k,n\in \mathbb N$.
\end{abstract}

\thanks{$^\dagger$ The second named author is affiliated with MaPhySto - A network in Mathematical Physics and Stochastics, which is funded by a grant from the Danish National Research Foundation.}

\maketitle

% \newpage
% \tableofcontents
% \newpage

\section{Introduction}
The quasi-nilpotent DT-operator $T$ was introduced by Dykema and the second author in \cite{DT}.
It can be described as the limit in $*$-moments for $n\to \infty$, of random matrices of the form
\begin{equation*}
  T^{(n)} = 
  \begin{pmatrix}
    0 & t_{1,2} & \cdots & t_{1,n} \\
    0 & \ddots & \ddots & \vdots \\
    \vdots & \ddots & \ddots & t_{n-1,n} \\
    0 & \cdots & 0 & 0
  \end{pmatrix}
\end{equation*}
where $\{\Re(t_{ij}),\Im(t_{ij})\}_{1\leq i < j\leq n}$ is a set of $n(n-1)$ independent identically distributed Gaussian random variables with mean $0$ and variance $\frac{1}{2n}$. More precisely, $T$ is an element in a finite von Neumann algebra, $M$, with a faithful normal tracial state, $\tau$, such that for all $s_1, s_2, \ldots, s_k\in \{1,*\}$, 
\begin{equation}
  \label{eq:1.1}
  \tau(T^{s_1} T^{s_2} \cdots T^{s_k}) = \lim_{n\to \infty} 
\mathbb E[\text{tr}_n((T^{(n)})^{s_1}(T^{(n)})^{s_2}\cdots 
(T^{(n)})^{s_k}  )], 
\end{equation}
where $\text{tr}_n$ is the normalized trace on $M_n(\mathbb C)$. Moreover the pair \\ $(T,W^*(T))$ is uniquely determined up to $*$-isomorphism by (\ref{eq:1.1}). The quasi-nilpotent DT-operator can be realized as an element in the free group factor, $L(\mathbb F_2)$, in the following way (cf. \cite[Sect. 4]{DT}): Let $(D_0, X)$ be a pair of free selfadjoint elements in a tracial $W^*$-probability space $(M,\tau)$, such that $\rm{d}\mu_{D_0}(t)=1_{[0,1]}(t)\rm{d}t$ and $X$ is semi-circular distributed, i.e. $\rm{d}\mu_X(t) = \frac{1}{2\pi} \sqrt{4-t^2}1_{[-2,2]}(t)\rm{d}t$.
Then $W^*(D_0,X) \simeq W^*(D_0)\star W^*(X)\simeq L(\mathbb F_2)$. 
Put
\begin{equation*}
  T_N=\sum_{j=1}^{2^N}p_{N,j}Xq_{N,j}
\end{equation*}
for $N=1,2,\ldots$, where 
\begin{equation*}
  p_{N,j}=1_{[\frac{j-1}{2^N},\frac{j}{2^N}]}(D_0), ~~~q_{N,j} = 1_{[\frac{j}{2^N},1]}(D_0),
\end{equation*}
for $j=1,2,\ldots, 2^N$. Then $(T_N)_{N=1}^\infty$ converges in norm to an operator $T\in W^*(D_0,X)$, and the $*$-moments of $T$ are given by (\ref{eq:1.1}), i.e. $T$ is a realization of the quasi-nilpotent DT-operator. In the notation of \cite[Sect. 4]{DT}, $T=\mathcal{UT}(X,\lambda)$, where $\lambda:L^\infty[0,1]\to W^*(D_0)$ is the $*$-isomorphism given by $\lambda(f)=f(D_0)$ for $f\in L^\infty([0,1])$. 
In the following we put $\mathcal D=W^*(D_0)\simeq L^\infty([0,1])$ and let $E_\mathcal D$ denote the trace-preserving conditional expectation of $W^*(D_0,X)$ onto $\mathcal D$.

In this paper we apply Voiculescu's $\mathcal R$-transform with amalgamation to compute various $*$-moments of $T$ and of operators closely related to $T$. First we compute in section 3 moments and the scalar valued $\mathcal R$-transform of $(T-\lambda 1)^*(T-\lambda 1)$ for $\lambda \in \mathbb C$. The specialized case of $\lambda =0$ was treated in \cite{DT} by more complicated methods. In section 4 we consider the operator
\begin{equation*}
  T+\sqrt{\epsilon}Y,
\end{equation*}
where $Y$ is a circular operator $*$-free from $T$ and $\epsilon>0$. By random matrix considerations it is easily seen, that if $T_1$ and $T_2$ are two quasi-nilpotent DT-operators, which are $*$-free with respect to amalgamation over the same diagonal, $\mathcal D$, then $T+\sqrt{\epsilon}Y$ has the same $*$-distribution as $S=\sqrt{a}T_1+\sqrt{b}T_2$, when $a=1+\epsilon$ and $b=\epsilon$ (cf. \cite{LA}).
We use this fact to prove, that the Brown measure of $T+\sqrt{\epsilon}Y$ is equal to the uniform distribution on the closed disc 
$\overline{B}(0,\log(1+\frac{1}{\epsilon})^{-\frac{1}{2}})$ in the complex plane. Moreover we show, that the spectrum of $T+\sqrt{\epsilon}Y$ is equal to this disc, and that $T+\sqrt{\epsilon}Y$ is not a DT-operator for any $\epsilon>0$. 

In \cite{DT} it was conjectured, that
\begin{equation}
  \label{eq:1.2}
  \tau(((T^*)^k T^k)^n )=\frac{n^{nk}}{(nk+1)!}
\end{equation}
for $n,k\in\mathbb N$. This formula was proved by \'Sniady in \cite{Sniady}. \'Sniady's proof of (\ref{eq:1.2}) is based on Speicher's combinatorial approach to free probability with amalgamation from \cite{AMS}. The key step in the proof of (\ref{eq:1.2}) was to establish a recursion formula for the $\mathcal D$-valued moments,
\begin{equation}
  \label{eq:1.3}
  E_\mathcal D\left(((T^*)^k T^k)^n  \right)
\end{equation}
for each fixed $k\in\mathbb N$. \'Sniady's recursion formula for the $\mathcal D$-valued moments (\ref{eq:1.3}), was later used by Dykema and the second author to prove, that
\begin{equation*}
  W^*(T) = W^*(D_0,X) \simeq L(\mathbb F_2)
\end{equation*}
and that $T$ admits a one parameter family of non-trivial hyperinvariant subspaces (cf. \cite{invsub}). In section 5 and section 6 of this paper we give a new proof of \'Sniady's recursion formula for the $\mathcal D$-valued moments (\ref{eq:1.3}), which at the same time gives a new proof of (\ref{eq:1.2}). The new proof is based on Voiculescu's $\mathcal R$-transform with respect to amalgamation over $M_{2k}(\mathcal D)$, the algebra of $2k\times 2k$ matrices over $\mathcal D$.

\section{Preliminaries}

In this section we give a few preliminaries on amalgamated probability
theory. Let $\mathcal A$ be a unital Banach algebra, and let $\mathcal B$ be a
Banach-sub-algebra containing the unit of $\mathcal A$. Then a map, $E_\mathcal B:\mathcal A\to
\mathcal B$, is a conditional expectation if 
\begin{enumerate}[(a)]
\item $E_\mathcal B$ is linear,
\item $E_\mathcal B$ preserves the unit i.e. $E_\mathcal B(1)=1$
\item and $E_\mathcal B$ has the $\mathcal B$, $\mathcal B$ bi-module property
i.e. $E_\mathcal B(b_1a b_2)=b_1
ab_2$ for all $b_1,b_2\in \mathcal B$ and $a\in \mathcal A$. 
\end{enumerate}

If $\mathcal B$, $\mathcal A$ and $E_\mathcal B$ are as above we say
that $(\mathcal B\subset \mathcal A, E_\mathcal B)$ is a $\mathcal
B$-probability space. If $\phi:\mathcal A\to \mathbb C$ is a state on
$\mathcal A$ which respects $E_\mathcal B$, i.e. $\tau= \tau\circ
E_\mathcal B$, we say that $(\mathcal B\subset \mathcal A, E_\mathcal B)$
is compatible to the (non-amalgamated) free probability space $(\mathcal 
A,\phi)$.

If $(\mathcal B\subset A, E_\mathcal B)$ is a $\mathcal B$-probability
space and $a\in \mathcal A$ is a fixed variable, we define the
amalgamated Cauchy transform of $a$ by
\begin{equation*}
G_a(b) = E_\mathcal B((b-a)^{-1}).
\end{equation*}
for $b\in \mathcal B$ and $b-a\in \mathcal B_{\textrm{inv}}$. The Cauchy transform is 1-1 in \\ $\{b\in \mathcal B_{\text{inv}}| \norm{b^{-1}} < \epsilon\}$ for $\epsilon$ sufficiently small and Voiculescu's amalgamated $\mathcal R$-transform \cite{Voi} is now defined for $a\in \mathcal A$ by
\begin{equation}
\mathcal R_a(b) = G^{\langle -1 \rangle}_a(b) - b^{-1},
\end{equation}
for $b$ being an invertible element of $\mathcal B$ suitably close to
zero. It turns out that this definition coincides on invertible
element with Speicher's definition of the amalgamated $\mathcal
R$-transform (cf. \cite[Th. 4.1.2]{AMS} and \cite{RS});
  \begin{equation}
    \label{R-trans}
    \mathcal R_a(b) = \sum_{n=1}^\infty \kappa_n^\mathcal B
    (a\otimes_\mathcal B ba\otimes_\mathcal B \cdots \otimes_\mathcal B 
    ba). 
  \end{equation}

We will need the following useful lemma for solving equations involving the amalgamated $R$-transform and Cauchy-transform.
\begin{lemma} \label{RCtheorem}
  Let $(\mathcal B \subset \mathcal A,E_\mathcal B)$ be a $\mathcal B$-probability space, and let $a\in \mathcal A$. Then there exists $\delta>0$ such that if $b\in \mathcal B$ is invertible, $\norm{b}< \delta$, $|\mu|>\frac{1}{\delta}$ and
  \begin{equation*}
    \mathcal R_a^\mathcal B(b) + b^{-1} = \mu 1_\mathcal A
  \end{equation*}
then $b = G_a^\mathcal B (\mu 1_\mathcal A)$.
\end{lemma}

\begin{proof}
  Let $\delta=\frac{1}{11\norm{a}}$ and define $g_b(b)= G_a^\mathcal B(b^{-1})$. By \cite[Prop. 2.3]{RS} we know that $g_a$ maps $\mathcal B(0,\frac{1}{4\norm{a}})$ bijectively onto a neighboorhood of zero containing $\mathcal B(0,\frac{1}{11\norm{a}})$ and furthermore that
  \begin{equation*}
    g_a^{\langle -1 \rangle}\left(\mathcal B(0,\tfrac{1}{11\norm{a}})_{\text{inv}}\right) \subseteq \mathcal B(0,\tfrac{2}{11\norm{a}})_{\text{inv}}.
  \end{equation*}
By definition we know that
\begin{equation*}
\mathcal R_a^\mathcal B(b) = {G_a^{\mathcal B}}^{\langle-1\rangle}(b) + b^{-1} = \left(g_a^{\langle -1 \rangle}(b)\right)^{-1} + b^{-1}
\end{equation*}
so if $\mathcal R_a(b) + b^{-1} = \mu 1_\mathcal A$ then
\begin{equation*}
  \mu 1_\mathcal A = g^{\langle -1\rangle }(b) - b^{-1} + b^{-1}  =\left(g_a^{\langle -1\rangle }(b)\right)^{-1}
\end{equation*}
and thus
\begin{equation}\label{galigning}
  g^{\langle -1 \rangle}_a(b) = \tfrac{1}{\mu} 1_\mathcal A.
\end{equation}
If $|\mu|> \frac{1}{\delta}$ then especially $\frac{1}{|\mu|}< \frac{1}{4\norm{a}}$ so $\tfrac{1}{\mu}1_\mathcal A$ is in the bijective domain of $g_a$, so applying $g_a$ on both sides of (\ref{galigning}) we get exactly
\begin{equation*}
  G^{\mathcal B}_a(\mu 1_\mathcal A) = g_a(\tfrac{1}{\mu} 1_\mathcal A) = b
\end{equation*}
since also $\norm{b} < \frac{1}{11\norm{a}}$. 
\end{proof}

If $a\in \mathcal A$ is a random variable in the $\mathcal
B$-probability space $(\mathcal B\subset\mathcal A,E_\mathcal B)$, then
following Speicher we define $a$ to be $\mathcal B$-Gaussian \cite[Def
4.2.3]{AMS} 
if only $\mathcal 
B$-cumulants of length 2 survive. From (\ref{R-trans}) it follows that in
this case the $\mathcal R$-transform has a particularly simple form,
namely, 
\begin{equation} \label{R-D-Gaussian}
\mathcal R_a(b)= \kappa_2^\mathcal B(a\otimes_\mathcal B ba)= E_\mathcal B(aba).
\end{equation}

In the following theorem (which is probably not a new one we just could
not find a proper reference) concerning cumulants we have adopted the
notation of Speicher from \cite{AMS}. 
\begin{lemma} \label{tensorcumulant}
  Let $N\in \mathbb N$ and let $(\mathcal B \subset \mathcal
  A,E_\mathcal B)$ be a $\cd$-probability space. Then $(M_N(\mathcal
  B)\subset M_N(\mathcal 
  A),E_{M_n(\mathcal B)})$ is a
  $M_N(\mathcal B)$-probability space with cumulants determined by the following
  formula: \\ 
\begin{multline*} 
\kappa_n^{M_N(\mathcal B)}((m_1\otimes
a_1)\otimes_{M_N(\mathcal B)}\cdots  \otimes_{M_N(\mathcal B)} (m_n\otimes a_n))   \\ =
(m_1 \cdots m_n) \otimes \kappa^\mathcal B_n (a_1
\otimes_\mathcal B \cdots \otimes_\mathcal B a_n) 
\end{multline*} 
when $m_1,\ldots ,m_n 
\in M_N(\mathbb C)$ and $ a_1 ,\ldots , a_n \in \mathcal A$.
\end{lemma}
We have of course made the identification $M_N(\mathcal A)\cong
M_N(\mathbb C)\otimes \mathcal A$.

\begin{proof}
  Since $M_N(\mathbb C)\subset M_N(\mathcal B)$ we observe that \\
\begin{multline*} 
\kappa_n^{M_N(\mathcal B)}((m_1\otimes
a_1)\otimes_{M_N(\mathcal B)}\cdots \otimes_{M_N(\mathcal B)}
(m_n\otimes a_n)) \\ =
((m_1 \cdots m_n)  \otimes  1) \cdot \\ 
\kappa^{M_N(\mathcal B)}_n ((1\otimes a_1) \otimes_{M_N(\mathcal B)} \cdots
\otimes_{M_N(\mathcal B)} (1\otimes a_n)).
\end{multline*}
To finish the proof we claim that
\begin{multline} \label{kappa}
  \kappa^{M_N(\mathcal B)}_n ((1\otimes a_1) \otimes_{M_N(\mathcal
  B)} \cdots \otimes_{M_N(\mathcal B)}
  (1\otimes a_n)) = \\ 1\otimes \kappa^\mathcal B_n (a_1\otimes_\mathcal B \cdots \otimes_\mathcal B a_n).
\end{multline}
The case $n=1$ is obvious since 
\begin{equation*}
1_N\otimes \kappa_1^\mathcal B(a_1) = 1_N \otimes E_\mathcal B(a_1) =
E_{M_N(\mathcal B)}(1\otimes a_1)  = \kappa_1^{M_N(\mathcal
B)}(1\otimes a_1).  
\end{equation*} 
Now assume that the claim is true for $1,2,\ldots,n-1$. Then
$(\ref{kappa})$ has an obvious extension to noncrossing partions of
length less than or equal to $n-1$. Hence     
\begin{multline*}
1_N\otimes \kappa_n^\mathcal B(a_1\otimes_\mathcal B \cdots \otimes_\mathcal B
 a_n)   \\ =
1_N\otimes E_\cd(a_1\cdots a_n)   - \sum_{\pi\in NC(n), \pi\neq
1_n} 1\otimes \kappa_\pi^\cd (a_1\otimes_\cd \cdots
\otimes_\cd a_n) \\ 
 = E_{M_N(\cd)}(1\otimes_{M_N(\mathcal B)}a_1\cdots a_n)  \\ -
 \sum_{\pi\in NC(n),\pi \neq 1_n} 
\kappa_\pi^{M_N(\cd)}((1\otimes 
a_1)\otimes_{M_N(\cd)}\cdots\otimes_{M_N(\cd)} (1\otimes a_n))  \\ 
= \kappa_n^{M_N(\cd)}((1\otimes a_1)\otimes_{M_N(\cd)} \cdots 
\otimes_{M_N(\cd)} (1\otimes a_n)). 
\end{multline*}
By induction this proves the lemma. 
\end{proof}

Assume that $\mathcal M$ contains a pair $(D_0, X)$ of $\tau$-free selfadjoint elements such that $\rm{d}\mu_{D_0}(t)=1_{[0,1]}(t)\rm{d}t$ and $X$ is a semicircular distributed. Put $\mathcal D=W^*(D_0)$. Then $\lambda:L^\infty([0,1])\to \mathcal D$ given by
\begin{equation*}
  \lambda(f)= f(D_0),
\end{equation*}
for $f\in L^\infty([0,1])$ is a $*$-isomorphism of $L^\infty([0,1])$ onto $\mathcal D$ and
\begin{equation*}
  \tau \circ \lambda(f) = \int_0^1 f(t)\rm{d}t, ~~f\in L^\infty([0,1]).
\end{equation*}
We will identify $\mathcal D$ with $L^\infty([0,1])$ and thus consider elements of $\mathcal D$ as functions. As explained in the introduction, we can realize the quasi-nilpotent DT-operator as the operator $T=\mathcal{UT}(X, \lambda)$ in $W^*(D_0, X)\simeq L(\mathbb F_2)$.

Define for $f\in \mathcal D\simeq L^\infty([0,1])$
\begin{eqnarray}
(L^*(f))(x) := \int_0^x f(t)\text{d}t & \textrm{and} & (L(f))(x) := \int_x^1 f(t)\text{d}t.
\end{eqnarray}
From the appendix of \cite{invsub} it follows that $(T,T^*)$ is a $\mathcal D$-Gaussian pair and that the covariances of
$(T,T^*)$ are given by the following lemma
\begin{lemma}\cite[Appendix]{invsub}
Let $f\in \mathcal D$. Then
\begin{eqnarray*}
E_\mathcal D (T f T^*) =  L(f) &  \textrm{and} & E_\mathcal D (T^*
f T)= L^*(f) 
\end{eqnarray*}
and $E_\mathcal D(TfT)=E_\mathcal D(T^*fT^*)=0$.
\end{lemma}

\section{Moments and $\mathcal R$-transform of $(T-\lambda 1)^*(T-\lambda 1)$} 
 \label{Sniadyresult}

Let $T$ be the quasi-diagonal DT-operator and define
\begin{equation*}
  \tilde{T} = 
  \begin{pmatrix}
    0 & T^* \\ T & 0
  \end{pmatrix}.
\end{equation*}
Since $(T, T^*)$ is a $\mathcal D$-Gaussian pair, it follows from lemma \ref{tensorcumulant}, that cumulants of the form 
\begin{equation*}
  \kappa_n^{M_2(\mathcal D)}((m_1\otimes a_1)\otimes_{M_2(\mathcal D)} \cdots \otimes_{M_2(\mathcal D)}(m_n\otimes a_n) )
\end{equation*}
vanishes when $n\neq 2$, $m_1, m_2, \ldots, m_n\in M_2(\mathbb C)$ and $a_1, a_2,\ldots , a_n \in \{T, T^*\}$.
Hence by the linearity of $\kappa_n^{M_2(\mathcal D)}$,
\begin{equation*}
  \kappa_n^{M_2(\mathcal D)}(\tilde{T}\otimes_{M_2(\mathcal D)} \tilde{T}\otimes_{M_2(\mathcal D)} \cdots \otimes_{M_2(\mathcal D)}\tilde{T})=0  
\end{equation*}
when $n\neq 2$, i.e. $\tilde{T}$ is a $M_2(\mathcal D)$-Gaussian element in $M_2(\mathcal M)$ under the conditional expectation $E_{M_2(\mathcal D)}:M_2(\mathcal M) \to M_2(\mathcal D)$ given by 
\begin{equation*}
  E_{M_2(\mathcal D)}:\begin{pmatrix}
    a_{11} & a_{12} \\
    a_{21} & a_{22}
  \end{pmatrix}\mapsto 
  \begin{pmatrix}
    E_\mathcal D(a_{11}) & E_\mathcal D(a_{12}) \\
    E_\mathcal D(a_{21}) & E_\mathcal D(a_{22})
  \end{pmatrix}.
\end{equation*}
Since $\tilde{T}$ is $M_2(\mathcal D)$-Gaussian the $\mathcal R$-transform of $\tilde{T}$ is by (\ref{R-D-Gaussian}) the linear mapping $M_2(\mathcal D) \to M_2(\mathcal D)$ given by
\begin{multline*}
  \mathcal R^{M_2(\mathcal D)}_{\tilde{T}}(z) = E_{M_2(\mathcal D)}(\tilde{T}z\tilde{T}) 
\\ = E_{M_2(\mathcal D)}\left(
  \begin{pmatrix}
   0 & T^* \\
   T & 0
  \end{pmatrix}
  \begin{pmatrix}
    z_{11} & z_{12} \\
    z_{21} & z_{22}
  \end{pmatrix}
  \begin{pmatrix}
    0 & T^* \\
   T & 0 
  \end{pmatrix}
\right) \\
= E_{M_2(\mathcal D)}\left(
  \begin{pmatrix}
  T^*z_{22}T & 0 \\
  0 & Tz_{11}T^* 
  \end{pmatrix}
\right) \\ 
= 
\begin{pmatrix}
  E_\mathcal D(T^* z_{22} T) & 0 \\
  0  & E_\mathcal D(Tz_{11}T^*) 
\end{pmatrix} \\
= 
\begin{pmatrix}
  L^*(z_{22}) & 0 \\
  0 & L(z_{11})
\end{pmatrix}.
\end{multline*}
For $\lambda\in \mathbb C$, we put $T_\lambda-T_\lambda 1$ and define
\begin{equation*}
  \tilde{T}_\lambda = 
  \begin{pmatrix}
    0 & T^*_\lambda \\
    T_\lambda & 0
  \end{pmatrix}
= \tilde{T} - 
  \begin{pmatrix}
  0 & \overline{\lambda} 1 \\
 \lambda 1 & 0  
  \end{pmatrix}
\end{equation*}
Since $ \begin{pmatrix}
  0 & \overline{\lambda} 1 \\
 \lambda 1 & 0  
  \end{pmatrix}\in M_2(\mathcal D)$ we have by $M_2(\mathcal D)$-freeness that the $\mathcal R$-transform is additive \cite[Th. 4.1.22]{AMS} i.e.
  \begin{equation*}
    \mathcal R_{\tilde{T}_\lambda}^{M_2(\mathcal D)}(z) = \mathcal R_{\tilde{T}}^{M_2(\mathcal D)} -   \begin{pmatrix}
    0 & \overline{\lambda} 1 \\
   \lambda 1 & 0
  \end{pmatrix} = 
  \begin{pmatrix}
    L^*(z_{22}) & -\overline{\lambda} 1 \\
   -\lambda 1 & L(z_{11})
  \end{pmatrix}.
\end{equation*}

One easily checks, that if $\delta\in \mathbb C$, $\delta\neq 0$, $\delta\neq -\frac{1}{|\lambda|^2}$ and $\mu\in\mathbb C$ is one of the two solutions to
\begin{equation*}
  \mu^2 = \frac{\text{e}^\sigma}{\sigma}(1+|\lambda|^2\sigma),
\end{equation*}
then
\begin{equation} \label{stjerneII}
\begin{cases}
  z_{11} = \mu\sigma \text{e}^{\sigma(x-1)} & \\
  z_{12} = -\overline{\lambda}\sigma & \\
  z_{21} = -\lambda\sigma & \\
  z_{22} = \mu\sigma \text{e}^{-\sigma x} & 
\end{cases}
\end{equation}
is a solution to
\begin{equation*}
\mathcal R_{\tilde{T}_\lambda}^{M_2(\mathcal D)}(z) + z^{-1} = \mu 1_2.
\end{equation*}
Here $x$ is the variable for the function in $\mathcal D=L^\infty([0,1])$.
In particular $z_{12}$ and $z_{21}$ are constant operators. If $\sigma \to 0$ then $|\mu|\to \infty$ and $\norm{z}\to 0$, so by lemma \ref{RCtheorem} there exists $\rho>0$ such that $|\sigma|<\rho$ implies
\begin{equation*}
G_{\tilde{T}_\lambda}^{M_2(\mathcal D)}(\mu 1_2) = 
\begin{pmatrix}
  z_{11} & z_{12} \\
 z_{21} & z_{22}
\end{pmatrix},
\end{equation*}
where $(z_{ij})_{i,j\in\{1,2\}}$ is given by (\ref{stjerneII}) and 
\begin{equation*}
  \mu = \pm \sqrt{\frac{\text{e}^\sigma}{\sigma}(1+|\lambda|^2\sigma)}.
\end{equation*}
On the other hand the Cauchy-transform of $\tilde{T}$ in $\mu 1_2$ is
\begin{multline*}
  \begin{pmatrix}
    z_{11} & z_{12} \\
    z_{21} & z_{22}
  \end{pmatrix}
= G_{\tilde{T}_\lambda}^{M_2(\mathcal D)}(\mu 1_2) \\
= E_{M_2(\mathcal D)}\left(
\left(\begin{pmatrix}
  \mu 1 & 0 \\
  0 & \mu 1
  \end{pmatrix} - 
  \begin{pmatrix}
    0 & T_\lambda^* \\
  T_\lambda & 0
  \end{pmatrix}\right)^{-1}
\right) \\
= E_{M_2(\mathcal D)}\left(
  \begin{pmatrix}
  \mu 1 & - T_\lambda^* \\
  -T_\lambda & \mu 1 
  \end{pmatrix}^{-1}
\right) \\
= E_{M_2(\mathcal D)}\left(
  \begin{pmatrix}
   \mu(\mu^2 1 - T^*_\lambda T_\lambda)^{-1} &  T^*_\lambda(\mu^2 1 - T_\lambda 
T^*_\lambda )^{-1} \\
  T_\lambda(\mu^2 1 - T^*_\lambda T_\lambda)^{-1} &  \mu(\mu^2 1 - T_\lambda T^*_\lambda)^{-1}
  \end{pmatrix} 
\right) .
\end{multline*} 
Thus
\begin{equation} \label{StjerneIII}
  \begin{cases}
    z_{11} = \mu E_\mathcal D ((\mu^2 1 - T_\lambda^*T_\lambda )^{-1}) & \\
    z_{12} = E_\mathcal D (T_\lambda^*(\mu^2 1 - T_\lambda T_\lambda^*)^{-1}) & \\
    z_{21} = E_\mathcal D (T_\lambda (\mu^2 1 - T_\lambda^*T_\lambda )^{-1}) & \\
    z_{22} = \mu E_\mathcal D ((\mu^2 1 - T_\lambda T_\lambda^*)^{-1}) &
  \end{cases}.
\end{equation}
Combining (\ref{stjerneII}) and (\ref{StjerneIII}) we have
\begin{equation} \label{Rcompu}
  \begin{cases}
E_\mathcal D ((\mu^2 1 - T_\lambda^*T_\lambda )^{-1}) = \sigma \text{e}^{\sigma(x-1)}  \\
E_\mathcal D (T_\lambda^*(\mu^2 1 - T_\lambda T_\lambda^*)^{-1}) = -\overline{\lambda}\sigma  \\
E_\mathcal D (T_\lambda (\mu^2 1 - T_\lambda^*T_\lambda )^{-1}) = -\lambda \sigma  \\
E_\mathcal D ((\mu^2 1 - T_\lambda T_\lambda^*)^{-1}) = \sigma \text{e}^{-\sigma x}
\end{cases}.
\end{equation}
We can now compute the $\mathcal R$-transform of $T_\lambda^* T_\lambda$ (wrt. $\mathbb C$) from (\ref{Rcompu}) and the defining equality for $\mu^2$.
\begin{multline*}
  \text{tr}\left(\left(\frac{\text{e}^\sigma}{\sigma}(1+|\lambda|^2\sigma)1 - T_\lambda^* T_\lambda\right)^{-1}\right) = \int_0^1 \sigma \text{e}^{\sigma(x-1)}\text{d}x \\
 = \left[\text{e}^{\sigma(x-1)}\right]_0^1 = 1-\text{e}^{-\sigma}.
\end{multline*}
Thus 
\begin{equation*}
  G^{\mathbb C}_{T^*_\lambda T_\lambda}\left( \frac{\text{e}^\sigma}{\sigma}(1+|\lambda|^2\sigma) \right) = 1 - \text{e}^{-\sigma}
\end{equation*}
i.e. 
\begin{equation*}
  \mathcal R^{\mathbb C}_{T^*_\lambda T_\lambda}(1-\text{e}^{-\sigma}) = \frac{\text{e}^\sigma}{\sigma}(1+|\lambda|^2\sigma) - \frac{1}{1-\text{e}^{-\sigma}}
\end{equation*}
for $\sigma$ in a neighboorhood of zero. Substituting $z=1-\text{e}^{-\sigma}$ we get $\sigma = -\log(1-z)$, so
\begin{equation*}
  \mathcal R^{\mathbb C}_{T^*_\lambda T_\lambda}(z) = - \frac{1}{(1-z)\log(1-z)}(1-|\lambda|^2\log(1-z)) -\frac{1}{z}. 
\end{equation*}
Hence we have proved the following extension of \cite[Theorem 8.7(b)]{DT}:
\begin{theorem} \label{Theorem:3.1a}
  Let $T$ be the quasinilpotent \rm{DT}-operator. Let $\lambda\in \mathbb C$ and put $T_\lambda = T-\lambda 1$. Then
  \begin{equation*}
    \mathcal R^{\mathbb C}_{T^*_\lambda T_\lambda}(z) = -\frac{1}{(1-z)\log(1-z)} -\frac{1}{z} 
+ \frac{|\lambda|^2}{1-z}
  \end{equation*}
for $z$ in some neighborhood of $0$.
\end{theorem}

We next determine the $\mathcal D$-valued (resp. $\mathbb C$-valued) moments of $T_\lambda^* T_\lambda$ for all $\lambda \in \mathbb C$. The special case $\lambda =0$ was treated in \cite[Theorem 5]{Sniady} (resp. \cite[Theorem 8.7(a)]{DT}) by different methods.

\begin{theorem} \label{Pr3.1}
  Let $\lambda \in \mathbb C$ and let $T, T_\lambda$ be as in theorem \ref{Theorem:3.1a}
  \begin{enumerate}[(a)]
  \item Let $Q_n$ be the sequence of polynomials on $\mathbb R$ uniquely determined by the following recursion formula
 \begin{equation} \label{3.4a}
    \begin{cases}
      Q_0(x)=1, &  \\
      Q_{n+1}(x) = |\lambda|^2Q_n(x+1) + \int_0^x Q_n(y+1)\text{d}y & \text{for }n\geq 1.
    \end{cases}
  \end{equation}
Then 
\begin{equation*}
  E_\mathcal D((T^*_\lambda T_\lambda)^n)(x)=Q_n(x),~~x\in [0,1],~~n\in \mathbb N.
\end{equation*}
\item 
  \begin{equation*}
    \tau((T^*_\lambda T_\lambda)^n)= \sum_{k=0}^n \frac{n^k}{(k+1)!}\binom{n}{k} |\lambda|^{2n-2k},~~n\in \mathbb N.
  \end{equation*}
  \end{enumerate}
\end{theorem}

\begin{proof}
By (\ref{Rcompu}), we have
\begin{equation}
  \label{eq:3.5a}
  E_\mathcal D\bigl((\tfrac{\rm{e}^\sigma}{\sigma}
(1+|\lambda|^2\sigma)1-T^*_\lambda T_\lambda)^{-1}\bigr)= \sigma \text{e}^{\sigma(x-1)}
\end{equation}
for $\sigma\in B(0,\rho)\setminus\{0\}$ for some $\rho>0$. Put
\begin{equation*}
  \psi(\sigma) = \frac{\sigma}{\rm{e}^\sigma (1+|\lambda|^2 \sigma)}, ~~ \sigma\in \mathbb C\setminus \{-\tfrac{1}{|\lambda|^2}\}.
\end{equation*}
Since $\psi(0)=0$ and $\psi'(0)=1$, $\psi$ has an analytic invers $\psi^{\langle -1 \rangle}$ defined in a neighborhood $B(0,\delta)$ of $0$, and we can choose $\delta>0$, such that $\psi^{\langle -1\rangle}(B(0,\delta))\subset B(0,\rho)$. By (\ref{eq:3.5a})
\begin{equation*}
  E_\mathcal D((\tfrac{1}{t} 1 - T^*_\lambda T_\lambda)^{-1})= \psi^{\langle -1 \rangle}(t)\text{e}^{\psi^{\langle -1 \rangle}(t)(x-1)}
\end{equation*}
for $t\in B(0,\delta)\setminus \{0\}$. By power series expansion of the left hand side, we get
\begin{equation}
  \label{eq:3.6a}
  \sum_{n=0}^\infty t^{n+1} E_\mathcal D((T^*_\lambda T_\lambda)^n)=\psi^{\langle -1 \rangle}(t)\text{e}^{\psi^{\langle -1 \rangle}(t)(x-1)}
\end{equation}
for $t\in B(0,\delta')$, where $0<\delta'\leq\delta $ and where the LHS of (\ref{eq:3.6a}) is absolutely convergent in the Banach space $L^\infty([0,1])$. Hence by Cauchy's integral formulas
\begin{equation}\label{eq:3.7a}
  E_\mathcal D((T^*_\lambda T_\lambda)^n)=\frac{1}{2\pi\rm{i}}\int_C \frac{\psi^{\langle -1 \rangle}(t)\text{e}^{\psi^{\langle -1 \rangle}(t)(x-1)}}{t^{n+2}}\rm{d}t
\end{equation}
as a Banach space integral in $L^\infty([0,1])$, where $C= \partial B(0,r)$ with positive orientation and $0<r<\delta'$.
For each fixed $x\in \mathbb R$
\begin{equation*}
  t\mapsto \psi^{\langle -1 \rangle}(t)\text{e}^{\psi^{\langle -1 \rangle}(t)(x-1)}
\end{equation*}
is an analytic function in $B(0,\delta')$ which is $0$ for $t=0$. Hence the function has a power series expansion of the form
\begin{equation}\label{eq:3.8a}
  \psi^{\langle -1 \rangle}(t)\text{e}^{\psi^{\langle -1 \rangle}(t)(x-1)}=\sum_{n=0}^\infty Q_n(x)t^{n+1}\end{equation}
for $t\in B(0,\delta')$, where the numbers $(Q_n(x))_{n=0}^\infty$ are given by
\begin{equation}
  \label{eq:3.9a}
  Q_n(x) = \frac{1}{2\pi\rm{i}}\int_C \frac{\psi^{\langle -1 \rangle}(t)\text{e}^{\psi^{\langle -1 \rangle}(t)(x-1)}}{t^{n+2}}\rm{d}t.
\end{equation}
In particular the $Q_n$'s are continuous functions of $x\in \mathbb R$. Substituting $\sigma=\psi(t)$ in (\ref{eq:3.8a}) we get
\begin{equation*}
  \sum_{n=0}^\infty Q_n(x)\psi(\sigma)^{n+1}= \sigma \text{e}^{\sigma(x-1)}
\end{equation*}
for $\sigma \in B(0,\rho')$, where $\rho'\in (0,\rho)$. Put
\begin{equation*}
  \begin{cases}
    R_0(x)=0 \\
    R_{n+1}(x) = |\lambda|^2 Q_n(x+1)+ \int_0^x Q_n(x)\rm{d}y, & n\geq 0.
  \end{cases}
\end{equation*}
Then
\begin{multline*}
  \sum_{n=0}^\infty R_n(x) \psi(\sigma)^{n+1} =
\psi(\sigma) \left(1+ \sum_{n=0}^\infty R_{n+1}(x)\psi(\sigma)^{n+1}\right) \\
= \psi(\sigma)\left(1+|\lambda|^2 \left(\sum_{n+0}^\infty Q_n(x+1)\right) +\int_0^x 
\left( \sum_{n+0}^\infty Q_n(y+1)\right)\rm{d}y\right) \\
=\psi(\sigma)\left(1+|\lambda|^2 \sigma \text{e}^{\sigma x}+ \int_0^x \sigma \text{e}^{\sigma y}\rm{d}y \right) ~~~~~~~~~~~~~~~~~~~~~~~~~~~~~~~~~~~\\
= \psi(\sigma)(|\lambda|^2\sigma +1) \text{e}^{\sigma x} 
= \sigma \text{e}^{\sigma(x-1)} 
= \sum_{n=0}^\infty Q_n(x) \psi(\sigma)^{n+1}~~~~~~~~~~~~~~~~
\end{multline*}
for all $\sigma\in B(0,\rho')$. Since $\psi(B(0, \rho'))$ is an open neighborhood of $0$ in $\mathbb C$, it follows that $R_n(x)=Q_n(x)$ for all $n\in \mathbb N$ and all $x\in \mathbb R$.

Hence $(Q_n(x))_{n=0}^\infty$ is the sequence of polynomials given by the recursive formula (\ref{3.4a}). Moreover by (\ref{eq:3.7a}) and (\ref{eq:3.9a}), $E_\mathcal D((T^*_\lambda T_\lambda)^n)=Q_n$ as functions in $L^\infty([0,1])$. This proves (a).

(b) By (\ref{eq:3.7a}), we have
\begin{equation*}
  \tau((T^*_\lambda T_\lambda)^n) = \int_0^1 E_\mathcal D((T^*_\lambda T_\lambda)^n)\rm{d}x
= \frac{1}{2\pi\rm{i}} \int_C \frac{1-\text{e}^{-\psi^{\langle -1 \rangle}(t)}}{t^{n+2}}\rm{d}t.
\end{equation*}
Note that $C'=\psi(C)$ is a positively oriented simple path around $0$. Hence by the substitution $t= \psi(\sigma)$, we get
\begin{eqnarray*}
  \tau((T^*_\lambda T_\lambda)^n) &=& \frac{1}{2\pi\text{i}}
\int_{C'} \frac{\psi'(\sigma)}{\psi(\sigma)^{n+2}}(1-\text{e}^{-\sigma}) \text{d}\sigma \\ 
&=& \frac{1}{2\pi\text{i}}
\int_{C'} \frac{1}{n+1}\frac{1}{\psi(\sigma)^{n+1}}\frac{\text{d}}{\text{d}\sigma}(1-\text{e}^{-\sigma}) \text{d}\sigma \\
&=& \frac{1}{2\pi\text{i}(n+1)}
\int_{C'} \frac{1}{\psi(\sigma)^{n+1}}\text{e}^{-\sigma} \text{d}\sigma \\
&=&\frac{1}{n+1}\left(\frac{1}{2\pi\text{i}}
\int_{C'} \frac{\text{e}^{n\sigma}(1+|\lambda|^2\sigma)^{n+1}}{\sigma^{n+1}} \text{d}\sigma\right) \\
&=& \frac{1}{n+1} \text{Res}\left( \frac{\text{e}^{n\sigma}(1+|\lambda|^2\sigma)^{n+1}}{\sigma^{n+1}}   ,0\right)
\end{eqnarray*}
where the second equation is obtained by partial integration and the last equality is obtained by the Residue theorem.

The above Residue is equal to the coefficient of $\sigma^n$ in the Power series expansion of
\begin{equation*}
  \text{e}^{n\sigma}(1+|\lambda|^2\sigma)^{-1}=
\left(\sum_{k=0}^\infty \frac{(n\sigma)^k}{k!}\right)
\left(\sum_{i=1}^{n+1} \binom{n+i}{i}(|\lambda|^2\sigma )^i\right).
\end{equation*}
Hence
\begin{multline*}
  \tau((T^*_\lambda T_\lambda)^n)= \frac{1}{n+1} 
\sum_{k=0}^n \frac{n^k}{k!} \binom{n+1}{n-k} |\lambda |^{2(n-k)} \\
= \frac{1}{n+1} \sum_{k=0}^n \frac{n^k}{(k+1)!} \binom{n}{k}|\lambda|^{2n-2k}.
\end{multline*}
\end{proof}

\section{Spectrum and Brown-measure of $T+\sqrt{\epsilon}Y$}

Let $T$ be the quasinilpotent DT-operator and let $Y$ be a circular operator $*$-free from $T$. In this section we will show, that
\begin{equation*}
  \sigma(T+\sqrt{\epsilon}Y) = \overline{B}\Bigl(0, 
\frac{1}{\sqrt{\log(1+\epsilon^{-1})}}\Bigr)
\end{equation*}
and that the Brown-measure $\mu_{T+\sqrt{\epsilon}Y}$ is equal to the uniform distribution on $\overline{B}\left(0, 
\frac{1}{\sqrt{\log(1+\frac{1}{\epsilon})}}\right)$, i.e. it has constant density w.r.t. the Lebesque measure on this disk.
 
\begin{theorem}
  For every $\epsilon >0$
  \begin{equation}
  \sigma(T+\sqrt{\epsilon}Y) = \overline{B}\Bigl(0, 
\frac{1}{\sqrt{\log(1+\epsilon^{-1})}}\Bigr).
  \end{equation}
\end{theorem}

\begin{proof}
The result can be obtained by the method of Biane and Lehner 
\cite[Section 5]{BL}. Let $a\in \mathbb C\setminus\{0\}$. Since $\sigma(T)=\{0\}$ we can write 
\begin{equation*}
  a1-(T+\sqrt{\epsilon})Y = \sqrt{\epsilon}\bigl(\tfrac{1}{\sqrt{\epsilon}}1- 
Y(a1-T)^{-1}\bigr)(a1-T).
\end{equation*}
Hence 
\begin{equation} \label{spectruma}
  a\notin \sigma(T+\sqrt{\epsilon}Y) \text{ iff } \frac{1}{\sqrt{\epsilon}} \notin \sigma\bigl(Y(a1-T)^{-1} \bigr).
\end{equation}

Let $Y=UH$ be the polar decomposition of $Y$. Then $Y(a1-T)^{-1} = UH(a1-T)^{-1}$, where $U$ is $*$-free from $H(a1-T)^{-1}$. Hence $Y(a1-T)^{-1}$ is $R$-diagonal. Moreover, since $0\notin \sigma(Y)$, $Y(a1-T)^{-1}$ is not invertible, so by \cite[Prop. 4.6.(ii)]{HL}
\begin{equation} \label{spectrumb}
\sigma\bigl(Y(a1-T)^{-1} \bigr) = B\bigl(0, \norm{Y(a1-T)^{-1}}_2\bigr).
\end{equation}
By $*$-freeness of $Y$ and $(a1-T)^{-1}$ we have
\begin{multline} \label{orthoT}
  \norm{Y(a1-T)^{-1}}_2^2 = \norm{Y}_2^2 \norm{(a1-T)^{-1}}_2^2 \\
= \norm{(a1-T)^{-1}}_2^2 = \norm{\sum_{n=0}^\infty \frac{T^n}{a^{n+1}}}_2^2.
\end{multline}

Applying now \cite[lemma 7.2]{DT} to $D=1$ and $\lambda = \frac{1}{a}$ and $\mu=\delta_0$, we get
\begin{equation*}
  \norm{\sum_{n=0}^\infty \frac{T^n}{a^n}}_2^2 = |a|^2\left(
\exp{\left(\frac{1}{|a|^2}\right)-1}\right)
\end{equation*}

Hence by (\ref{orthoT})
\begin{equation*}
  \norm{Y(a1-T)^{-1}}_2^2 = \exp\bigl(\frac{1}{|a|^2}\bigr)-1.
\end{equation*}
Thus for $a\in \mathbb C\setminus \{0\}$ we get by (\ref{spectruma}) and (\ref{spectrumb}) 
\begin{multline*}
  a\notin \sigma(T+\sqrt{\epsilon}Y) \Leftrightarrow \frac{1}{\sqrt{\epsilon}}\notin \sigma\bigl(Y(a1-T)^{-1}\bigr) \\
\Leftrightarrow \frac{1}{\sqrt{\epsilon}} > \exp\bigl(\frac{1}{|a|^2}\bigr)-1\Leftrightarrow |a| > \frac{1}{\sqrt{\log(1+\frac{1}{\epsilon})}}.
\end{multline*}
Hence $\sigma(T+\sqrt{\epsilon}Y)\cup \{0\} = \overline{B}\left(0, \frac{1}{\sqrt{\log(1+\frac{1}{\epsilon})}}\right)$. Since $\sigma(T+\sqrt{\epsilon}Y)$ is closed it follows that $\sigma(T+\sqrt{\epsilon}Y) = \overline{B}\left(0, \frac{1}{\sqrt{\log(1+\frac{1}{\epsilon})}}\right)$.
\end{proof}

In order to compute the Brown measure of $T+\sqrt{\epsilon}Y$, we first observe that $T+\sqrt{\epsilon}Y$ has the same $*$-distribution as
\begin{equation*}
  S= \sqrt{a}T_1 + \sqrt{b}T_2^*
\end{equation*}
when $T_1$ and $T_2$ are two $\mathcal D$-free quasidiagonal operators 
and $a=1+\epsilon$ and $b=\epsilon$ \cite{LA}. We next compute the Brown measure of $S$ for all values of $a,b \in (0,\infty)$.

\begin{lemma} \label{Brownlemma}
  Let $\mu_Q$ be the Brown measure of an operator $Q$ in a tracial $W^*$-probability space $(M,\rm{tr})$. Let $r>0$ and assume that $\mu_Q(\partial B(0,r))=0$. Then 
  \begin{equation*}
    \mu_Q(B(0,r)) = -\frac{1}{2\pi} \lim_{\alpha\to 0^+}\Im \left(\int_{\partial B(0,r)} \rm{tr}((Q_\lambda^* Q_\lambda +\alpha1)^{-1}Q_\lambda^*) \rm{d}\lambda \right)
  \end{equation*}
where $Q_\lambda = Q-\lambda 1$ for $\lambda \in \mathbb C$.
\end{lemma}

\begin{proof}
  Let $\Delta:M\to [0,\infty)$ be the Fuglede-Kadison determinant on $M$, and put $L(\lambda) = \log \Delta(Q_\lambda)$ and
  \begin{equation*}
    L_\alpha(\lambda) = \log \Delta((Q_\lambda^*Q_\lambda + \alpha 1)^{1/2}) = \frac{1}{2}\text{tr}(\log(Q_\lambda^* Q_\lambda +\alpha 1))
  \end{equation*}
for $\lambda \in \mathbb C$.

Put $\lambda_1 = \Re \lambda, \lambda_2 = \Im \lambda$ and let $\nabla^2 = 
\frac{\partial^2}{\partial \lambda_1^2} +\frac{\partial^2}{\partial \lambda_2^2} $ denote the Laplace operator on $\mathbb C$. Then by \cite[Section 2]{La} $\nabla^2 L_\alpha \geq 0$ and for each $\alpha >0$, the measure 
\begin{equation} \label{mualphamaal}
  \mu_\alpha  =\frac{1}{2\pi}\nabla^2 L_\alpha(\lambda) \text{d}\lambda_1 \text{d}\lambda_2
\end{equation}
is a probability measure on $\mathbb C$. Moreover
\begin{equation}
  \label{weakmualpha}
  \lim_{\alpha \to 0} \mu_{\alpha} = \mu
\end{equation}
in the weak$^*$ topology on $\text{Prob}(\mathbb{C})$. Also from \cite[Section 2]{La} the gradient $(\frac{\partial}{\partial \lambda_1},\frac{\partial}{\partial \lambda_2})$ of $L_\alpha$ is given by
\begin{eqnarray} \label{tre1}
  \frac{\partial}{\partial \lambda_1}L_\alpha(\lambda) & = & 
       -\Re\bigl(\text{tr}(Q_\lambda(Q_\lambda^* Q_\lambda + \alpha 1)^{-1}\bigr)\\
\label{tre2} \frac{\partial}{\partial \lambda_2}L_\alpha(\lambda) & = & 
       -\Im\bigl(\text{tr}(Q_\lambda(Q_\lambda^* Q_\lambda + \alpha 1)^{-1}\bigr)\end{eqnarray}
By (\ref{weakmualpha})
\begin{equation*}
  \lim_{\alpha\to 0} \int_{\mathbb C} \phi \text{d}\mu_\alpha = \int_{\mathbb C} \phi \text{d}\mu
\end{equation*}
for all $\phi \in C_0(\mathbb C)$. Since $1_{B(0,r)}$ is the limit of an increasing sequence $(\phi_n)_{n=1}^\infty$ of $C_0(\mathbb C)$-functions with $0\leq \phi_n \leq 1$ for all $n\in \mathbb N$ it follows that
\begin{multline*}
  \mu_Q(B(0,r)) = \lim_{n\to \infty} \int_{\mathbb C} \phi_n \text{d}\mu_Q \\
= \lim_{n\to \infty} \left(\lim_{\alpha\to 0} \int_{\mathbb C} \phi_n\text{d}\mu_\alpha \right) \leq  \lim_{n\to \infty} \left(\liminf_{\alpha\to 0} \int_{\mathbb C} 1_{B(0,r)}\text{d}\mu_\alpha \right) \\ = \liminf_{\alpha\to 0} \mu_\alpha (B(0,r))
\end{multline*}
Writing $1_{\overline{B}(0,r)}$ as the limit of a decreasing sequence $(\psi_n)_{n=1}^\infty$ of $C_0(\mathbb C)$-functions, with $0\leq \psi_n \leq 1$, one gets in the same way
\begin{equation*}
  \mu_Q(\overline{B}(0,r)) \geq \limsup_{\alpha\to 0} \mu_\alpha(\overline{B}(0,r))
\end{equation*}
Hence if $\mu_Q(\partial B(0,r))=0$ we have
\begin{equation*}
  \limsup_{\alpha\to 0} \mu_\alpha(B(0,r)) \leq \mu_Q(B(0,r)) \leq  \liminf_{\alpha\to 0} \mu_\alpha(B(0,r)),
\end{equation*}
and therefore
\begin{equation*}
   \mu_Q(B(0,r)) =  \lim_{\alpha\to 0} \mu_\alpha(B(0,r)).
\end{equation*}

Using (\ref{mualphamaal}) together with Green's theorem applied to the vector-field $(P_\alpha,Q_\alpha)= (-\frac{\partial L_\alpha}{\partial \lambda_2},\frac{\partial L_\alpha}{\partial \lambda_1} )$ we get
\begin{eqnarray*}
  \mu_\alpha(B(0,r)) & = & \frac{1}{2\pi} \int_{B(0,r)}\nabla^2 L_\alpha(\lambda)\text{d}\lambda_1\text{d}\lambda_2 \\
& = &\frac{1}{2\pi} \int_{B(0,r)}\left(\frac{\partial Q_\alpha}{\partial \lambda_1} -\frac{\partial P_\alpha}{\partial \lambda_2} \right) \text{d}\lambda_1\text{d}\lambda_2 \\
&= &\frac{1}{2\pi} \int_{\partial B(0,r)}P_\alpha\text{d}\lambda_1+ Q_\alpha\text{d}\lambda_2  \\
&= &\frac{1}{2\pi} \int_{\partial B(0,r)}-\frac{\partial L_\alpha}{\partial \lambda_2}\text{d}\lambda_1+ \frac{\partial L_\alpha}{\partial \lambda_1}\text{d}\lambda_2 \\
& = & \Im \left(\frac{1}{2\pi}\int_{\partial B(0,r)} \left(\frac{\partial L_\alpha}{\partial \lambda_1}- \text{i} \frac{\partial L_\alpha}{\partial \lambda_2} \right)\left(       \text{d}\lambda_1 +\text{i}\text{d}\lambda_2\right)\right)
\end{eqnarray*}
By (\ref{tre1}) and (\ref{tre2}) 
\begin{equation*}
  \frac{\partial L_\alpha}{\partial \lambda_1}- \text{i} \frac{\partial L_\alpha}{\partial \lambda_2} 
= -\overline{\text{tr}(Q_\lambda(Q_\lambda^* Q_\lambda+\alpha 1)^{-1})}  
= -\text{tr}((Q_\lambda^* Q_\lambda+\alpha 1)^{-1}Q^*_\lambda).
\end{equation*}
Hence 
\begin{equation*}
  \mu_\alpha(B(0,r)) = -\Im\left(\frac{1}{2\pi}\int_{\partial B(0,r)} 
\text{tr}((Q^*_\lambda Q_\lambda+\alpha 1)^{-1}Q_\lambda^*)\text{d}\lambda\right)
\end{equation*}
which completes the proof of the lemma.
\end{proof}

Let $S=\sqrt{a}T_1 +\sqrt{b}T_2^*$ with $0<b<a$. Since $cS$ and $S$ have the same $*$-distribution for all $c\in \mathbb T$, the Brown measure $\mu_S$ of $S$ is rotation invariant (i.e. invariant under the transformation $z\mapsto cz,~z\in \mathbb C$ when $|c|=1$). Hence by lemma \ref{Brownlemma} we can compute $\mu_S$, if we can determine
\begin{equation*}
  \text{tr}((S_\lambda^*S_\lambda + \alpha 1)^{-1}S_\lambda^*)
\end{equation*}
for all $\lambda \in \mathbb C$, where $S_\lambda = S- \lambda 1$, and for all $\alpha$ in some interval of the form $(0,\alpha_0)$. This can be done by minor modifications of the methods used in section \ref{Sniadyresult}:

Put 
\begin{equation*}
  \tilde{S}_\lambda = 
  \begin{pmatrix}
    0 & S^*_\lambda \\
   S_\lambda & 0 
  \end{pmatrix}.
\end{equation*}
Then there exists a $\delta>0$ (depending on $a,b$ and $\gamma$) such that when $\norm{z}\leq \delta$ and $|\mu| > \frac{1}{\delta}$ the equality
\begin{equation} \label{Rdeltadef}
  \mathcal R_{\tilde{S}_\lambda}^{M_2(\mathcal D)}(z) + z^{-1} = \mu 1_2
\end{equation}
implies that
\begin{multline} \label{R2kryds2}
  z = G_{\tilde{S}_\lambda}^{M_2(\mathcal D)}(\mu 1_2) \\ 
= (\text{id}\otimes E_\mathcal D)
\begin{pmatrix}
  \mu(\mu^2 1 - S_\lambda^* S_\lambda)^{-1} & S_\lambda^*(\mu^2 1 - S_\lambda S_\lambda^*)^{-1}\\
   S_\lambda(\mu^2 1 - S_\lambda^* S_\lambda)^{-1}&  \mu(\mu^2 1 - S_\lambda S_\lambda^*)^{-1} 
\end{pmatrix}.
\end{multline}
Moreover, $\tilde{S} = 
\begin{pmatrix}
  0 & S^* \\ S & 0
\end{pmatrix}$
is $M_2(\mathcal D)$-Gaussian by lemma \ref{tensorcumulant} since $(T_1,T_2^*,T_2,T_2^*)$ is a $\mathcal D$-Gaussian set. Hence for $z= (z_{ij})_{i,j=1}^2 \in M_2(\mathcal D)$,
\begin{equation*}
  \mathcal R_{\tilde{S}}^{M_2(\mathcal D)}(z) = E_{M_2(\mathcal D)}(\tilde{S}z\tilde{S}) =   \begin{pmatrix}
    E_\mathcal D(S^*z_{22}S) & 0 \\
 0 & E_\mathcal D(Sz_{11}S^*)
  \end{pmatrix}.
\end{equation*}
Using that $(T_1,T_1^*)$ and $(T_2,T_2^*)$ have the same $\mathcal D$-distribution as $(T,T^*)$ and that $(T_1,T_1^*)$ and $(T_2,T_2^*)$ are two $\mathcal D$-free sets, we get
\begin{eqnarray*}
   E_\mathcal D(S^*z_{22}S) & = & (aL^*+bL)(z_{22})\\
   E_\mathcal D(Sz_{11}S^*) & = & (aL+bL^*)(z_{11}),
\end{eqnarray*}
where $L(f): x\mapsto \int_x^1 f(y)\text{d}y$ and $L^*(f):x\mapsto \int_0^x f(y)\text{d}y$ for $f\in \mathcal D$. 

Since $\tilde{S}_\lambda = \tilde{S} - 
\begin{pmatrix}
  0 & \overline{\lambda}1 \\
\lambda 1 & 0 
\end{pmatrix}$ it follows that 
\begin{equation*}
  \mathcal R_{\tilde{S}}^{M_2(\mathcal D)}(z) = 
  \begin{pmatrix}
    (aL+bL^*)(z_{22}) & -\overline{\lambda} 1 \\
  \lambda 1 & (aL^*+bL)(z_{11})
  \end{pmatrix}.
\end{equation*}
Thus (\ref{R2kryds2}) becomes
\begin{multline}
  \label{2*2ligninger}
    \begin{pmatrix}
    \mu 1 & 0 \\
    0 & \mu 1
  \end{pmatrix}
 \\ =  
\begin{pmatrix}
    (aL+bL^*)z_{22} & -\overline{\lambda}1 \\
  \lambda 1 & (aL^*+bL)(z_{11})
  \end{pmatrix} + \frac{1}{\text{det}(z)}
  \begin{pmatrix}
    z_{22} & -z_{12} \\
   -z_{21} & z_{11} 
  \end{pmatrix}.
\end{multline}

In analogy with section \ref{Sniadyresult}, we look for solutions $z_{ij} \in \mathcal D = L^\infty[0,1]$ of the form
\begin{equation} \label{Brownsol}
\begin{pmatrix}
  z_{11} & z_{12} \\
  z_{21} & z_{22}
\end{pmatrix}=
\begin{pmatrix}
  c_{11}\exp(\sigma x) & c_{12} \\
  c_{21} & c_{22}\exp(-\sigma x)
\end{pmatrix},
\end{equation}
where $\sigma\in \mathbb C$ and $c=
\begin{pmatrix}
  c_{11} & c_{12} \\
  c_{21} & c_{22} 
\end{pmatrix}\in \text{GL}(2,\mathbb C)$.
It is easy to check that (\ref{Brownsol}) is a solution to (\ref{2*2ligninger}) if the following 5 conditions are fulfilled:
\begin{eqnarray*}
  \det(c) & =& \frac{\sigma}{a-b} \\
  c_{11} & = & \frac{\sigma \mu}{ a\text{e}^\sigma -b} \\
  c_{12} & =& -\frac{\sigma\overline{\lambda}}{a-b} \\
  c_{21} & =& -\frac{\sigma\lambda}{a-b} \\
  c_{22} & = & \frac{\sigma \mu}{ a-b\text{e}^{-\sigma}} 
\end{eqnarray*}
The first of these conditions is consistent with the remaining 4 if and only if
\begin{equation*}
  \frac{(\sigma \mu)^2}{(a\text{e}^\sigma -b)(a-b\text{e}^{-\sigma})}- 
\frac{\sigma^2|\lambda|^2 }{(a-b)^2} = \frac{\sigma}{a-b}
\end{equation*}
which is equivalent to
\begin{equation}
  \label{mucondition}
  \mu^2 = \frac{(a\text{e}^\sigma -b)(a-b\text{e}^{-\sigma})(a-b+\sigma|\lambda|^2)}{\sigma(a-b)^2} .
\end{equation}
Put 
\begin{equation*}
  \sigma_0 := -\min\left\{\frac{a-b}{|\lambda|^2},\log\bigl(\frac{a}{b}\bigr) \right\}.
\end{equation*}
Then for $\sigma_0 < \sigma < 0$, the right hand side of (\ref{mucondition}) is negative. Let in this case $\mu(\sigma)$ denote the solution to (\ref{mucondition}) with positive imaginary part, i.e.
\begin{equation}
  \label{impartmu}
  \mu(\sigma) = \text{i}\frac{a\text{e}^{\sigma/2}-b\text{e}^{-\sigma/2}}{|\sigma|^{1/2}(a-b)}\sqrt{a-b+\sigma|\lambda|^2}
\end{equation}
for $\sigma_0<\sigma <0$. Then with
\begin{eqnarray*}
    c_{11}= \frac{\sigma \mu(\sigma)}{a\text{e}^\sigma -b} & ~~&
    c_{12}= -\frac{\sigma\overline{\lambda}}{a-b} \\
    c_{21}=-\frac{\sigma\lambda}{a-b} & ~~&
c_{22}= \frac{\sigma \mu(\sigma)}{a-b\text{e}^{-\sigma}}
\end{eqnarray*}
the matrix $z(\sigma) =
\begin{pmatrix}
  z_{11} & z_{12} \\
  z_{21} & z_{22}
\end{pmatrix}
$ given by (\ref{Brownsol}) is a solution to
\begin{equation*}
  \mathcal R_{\tilde{S}_\lambda}^{M_2(\mathcal D)}(z(\sigma)) + z(\sigma)^{-1} = \mu 1_2.
\end{equation*}
By (\ref{impartmu}) $\lim_{\sigma \to 0^-}|\mu(\sigma)|=\infty$ and 
$\lim_{\sigma\to 0^-} |\sigma \mu(\sigma)| =0$ and therefore $\lim_{\sigma\to 0^-} \norm{z(\sigma)} = 0$.

Hence for some $\sigma_1 \in (\sigma_0,0)$ we have $|\mu(\sigma)| > \frac{1}{\delta}$ and $\norm{z(\sigma)}>\delta$ when $\sigma \in (\sigma_1,0)$ where $\delta>0$ is the number described in connection with (\ref{Rdeltadef}). Thus 
\begin{equation}
  \label{6*}
  z(\sigma) = G_{\tilde{S}_\lambda}^{M_2(\mathcal D)}(\mu(\sigma)1_2)
\end{equation}
for $\sigma \in (\sigma_1,0)$. But since both $\sigma\mapsto z(\sigma)$ and $\sigma \mapsto \mu(\sigma)$ are analytic functions (of the real variable $\sigma$) it follows that (\ref{6*}) holds for all $\sigma\in (\sigma_0,0)$. Note that $\sigma\mapsto -\text{i}\mu(\sigma)$ is a continuous strictly positive function on $(\sigma_0,0)$, and 
\begin{eqnarray*}
  \lim_{\sigma\to 0^-}(-\text{i}\mu(\sigma)) = +\infty & ~~~ &  
\lim_{\sigma\to \sigma_0^+}(-\text{i}\mu(\sigma)) = 0.
\end{eqnarray*}
Hence for every fixed real number $\alpha>0$ we can chose $\sigma\in (\sigma_0,0)$, such that 
\begin{equation*}
 - \text{i}\mu(\sigma) = \sqrt{\alpha}.
\end{equation*}
Thus by (\ref{R2kryds2}) and (\ref{6*}) 
\begin{equation*}
  E_\mathcal D(S_\lambda^*(-\alpha 1- S_\lambda S_\lambda^*)^{-1}) = 
z(\sigma)_{12} = -\frac{\sigma\overline{\lambda}}{a-b}
\end{equation*}
which is a constant function in $L^\infty [0,1]$. Hence 
\begin{equation*}
  \text{tr}(S_\lambda^*(S_\lambda S_\lambda^*+\alpha 1)^{-1}) = 
 \frac{\sigma\overline{\lambda}}{a-b}
\end{equation*}
from which
\begin{equation*}
  \int_{\partial B(0,r)}\text{tr}( S_\lambda^*(S_\lambda S_\lambda^* +\alpha 1)^{-1})\text{d}\lambda = 2\pi\text{i}\frac{\sigma r^2}{a-b}
\end{equation*}
when $\sigma_0 < \sigma < 0$, where as before $\sigma_0 = -\min\left\{ \frac{a-b}{|\lambda|^2}, \log\bigl( \frac{a}{b}\bigr)\right\}$.

Now $\alpha\to 0^+$ corresponds to $\sigma \to \sigma_0^+$. Hence 
\begin{multline*}
\lim_{\alpha\to 0^+} \left( -\frac{1}{2\pi} \Im \int_{\partial B(0,r)} \text{tr}(S_\lambda^*(S_\lambda S_\lambda^* +\alpha 1)^{-1})\text{d}\lambda\right) \\ =
-\frac{\sigma_0 r^2}{a-b} = +\min\left\{1, r^2\frac{\log\bigl(\frac{a}{b}\bigr)}{a-b} \right\}.
\end{multline*}
Obeserve that $S^*_\lambda (S_\lambda S_\lambda^*+\alpha 1)^{-1} = (S_\lambda^* S_\lambda +\alpha 1)^{-1} S_\lambda^*$.
Thus by lemma \ref{Brownlemma} we have for all but countably many $r>0$, that 
\begin{equation*}
  \mu_S(B(0,r)) = \min\left\{ 1, r^2\frac{\log\bigl(\frac{a}{b}\bigr)}{a-b}\right\}  = 
  \begin{cases}
    r^2\frac{\log\bigl(\frac{a}{b}\bigr)}{a-b}, & r\leq \sqrt{\frac{a-b}{\log\bigl(\frac{a}{b}\bigr)}} \\
    1, &  r> \sqrt{\frac{a-b}{\log\bigl(\frac{a}{b}\bigr)}}
  \end{cases}.
\end{equation*}
Since the right hand side is a continuous function of $r$, the formula actually holds for all $r>0$. This together with the rotation invariance of $\mu_S$ shows, that $\mu_S$ is equal to the uniform distribution on 
$\overline{B}\left(0,\sqrt{\frac{a-b}{\log\bigl(\frac{a}{b}\bigr)}}\right)$, i.e. has constant density $\frac{1}{\pi}\frac{\log\bigl(\frac{a}{b}\bigr)}{a-b}$ on this ball, and vanishes outside the ball. Putting $a=1+\epsilon$ and $b=\epsilon$ we get in particular
\begin{theorem} \label{theorem:4.3a}
  The Brown measure of $T+\sqrt{\epsilon}Y$ is equal to the uniform distribution on $\overline{B}\Bigl(0, \frac{1}{\sqrt{\log(1+\epsilon^{-1})}}\Bigr)$.
\end{theorem}

The Brown mesure of $T+\sqrt{\epsilon}Y$ can be used to give an upper bound of the microstate entropy of $T+\sqrt{\epsilon}Y$. By \cite{SniadyDT} we have for $S\in \mathcal M$
\begin{equation} \label{SniadyDTeq}
  \chi(S) \leq \int_\mathbb C \int_\mathbb C \log|z_1-z_2| \text{d}\mu_S(z_1)\text{d}\mu_S(z_2) + \frac{5}{4} + \log(\pi\sqrt{2\text{od}_S})
\end{equation}
where $\mu_S$ is the Brown measure of $S$ on $\mathbb C$ and $\text{od}_S$ is the off-diagonality of $S$ defined by
\begin{equation} \label{offdiag1}
  \text{od}_S := \tau(SS^*) - \int_\mathbb C |z|^2\text{d}\mu_S(z). 
\end{equation}

\begin{lemma} \label{lemma:4.4n} For $R>0$ we have
  \begin{equation*}
    I:=\int_{B(0,R)}\int_{B(0,R)} \log|z_1-z_2|\text{d}z_1\text{d}z_2=\pi^2(R^2\log R-\tfrac{1}{4})
  \end{equation*}
\end{lemma}
\begin{proof}
  Polar substitution in $I$ gives
  \begin{equation*}
   I:= 4\pi^2\int_0^R \int_0^R \left(\frac{1}{2\pi}\int_0^{2\pi}\log|r-\text{e}^{\text{i}\theta}s|\text{d}\theta\right)r\text{d}rs\text{d}s.
  \end{equation*}
Let $0<s<r$. $z\mapsto \log|r-zs|$ is the real value of the complex holomorphic function $z\mapsto \text{Log}(r-zs)$, where $\text{Log}$ is the principal branch of the complex logarithm, so $z\mapsto \log|r-zs|$ is a harmonic function in $B(0,\frac{r}{s})$. By the mean value property of harmonic functions 
\begin{equation*}
  \frac{1}{2\pi}\int_0^{2\pi}\log|r-\text{e}^{\text{i}\theta}s|\text{d}\theta = \log(r),
\end{equation*}
so symmetry in $r$ and $s$ reduces $I$ to
\begin{multline*}
  I:=  4\pi^2\int_0^R \int_0^R \max\{\log(r),\log(s)\}r\text{d}rs\text{d}s \\
= 8\pi^2\int_0^R \left(\int_0^r \log(r)s\text{d}s\right)r\text{d}r \\
= 4\pi^2\int_0^R r^3\log(r)\text{d}r = \pi^2R^4(\log(R)-\tfrac{1}{4}).
\end{multline*}
\end{proof}

\begin{theorem} \label{theorem:4.5n}
  \begin{multline}
    \label{eq:4.18n}
    \chi(T+\sqrt{\epsilon}Y) \leq -\frac{1}{2}\log(\log(1+\epsilon^{-1}))-\frac{1}{4} + \log\pi \\ + \frac{1}{2}\log\left(1+2\epsilon -\frac{1}{\log(1+\epsilon^{-1})}\right).
  \end{multline}
\end{theorem}
\begin{proof}
  Let $\nu_R$ be the uniform distribution on $\overline{B}(0,R)$. Since $\nu_R$ has constant density $(\pi R^2)^{-1}$ on $\overline{B}(0, R)$, we have by lemma \ref{lemma:4.4n} 
  \begin{equation*}
    \int_\mathbb C \int_\mathbb C \log|z_1-z_2| \rm{d}\nu_R(z_1)\rm{d}\nu_R(z_2)= \log R - \frac{1}{4}.
  \end{equation*}
The Brown measure of $S=T+\sqrt{\epsilon}Y$ is $\mu_S=\nu_R$ with $R=\log(1+\epsilon^{-1})^{-\frac{1}{2}}$, and
\begin{equation*}
  \rm{od}_S = \frac{1}{2} +\epsilon - \int_\mathbb C |z|^2\rm{d}\nu_R  = \frac{1}{2} + \epsilon -\frac{R^2}{2}.
\end{equation*}
Hence by (\ref{SniadyDTeq}) 
\begin{equation*}
   \chi(T+\sqrt{\epsilon}Y) \leq \log R -\frac{1}{4} + \log\pi  
+ \frac{1}{2}\log(1+2\epsilon -R^2).
\end{equation*}
This proves (\ref{eq:4.18n}).
\end{proof}
In \cite{LA} the first author proved that the microstate-free analog, $\delta_0^*(T)$, of the free entropy dimension is equal to 2. From Theorem \ref{theorem:4.5n} one gets only the trivial estimate of the free entropy dimension $\delta_0(T)$, namely
\begin{equation}
  \label{eq:4.19n}
  \delta_0(T) \leq 2 +\lim_{\delta\to 0^+} \frac{\chi(T+\sqrt{2}\delta Y)}{|\log \delta|} =2.
\end{equation}
If $T+\sqrt{\epsilon}Y$ was a DT-operator for all $\epsilon>0$ then by \cite{SniadyDT} equality would hold in (\ref{eq:4.18n}), and hence also in (\ref{eq:4.19n}). In the rest of this section, we prove that unfortunately $T+\sqrt{\epsilon}Y$ is not a DT-operator for any $\epsilon>0$.

If $R=D+T$ is a $\text{DT}(\mu,1)$ operator it follows from 
\cite[lemma 7.2]{DT} that for $|\lambda| < \norm{R}^{-1}$,
\begin{equation*}
  \norm{\sum_{n=0}^\infty \lambda^n R^n}^2_2 = \frac{1}{|\lambda|^2}
\left(\text{exp}\left(\sum_{k,l=1}^\infty \lambda^{k+1}\overline{\lambda}^{l+1} 
M_\mu(k,l)-1\right)\right),
\end{equation*}
where $M_\mu(k,l) = \int_{\sigma(R)}z^k\overline{z}^l\text{d}\mu_R(z)$.

If thus $\mu_D$ is the uniform distribution on a disk with radius $d$ then
\begin{equation*}
  M_{\mu_D}(k,l) = 0
\end{equation*}
when $k\neq l$ and
\begin{multline*}
  M_{\mu_D}(k,k) = \frac{1}{\pi d^2}\int_{B(0,d)}|z|^{2k} \text{d}z \\
 = \frac{2\pi}{\pi d^2}\int_0^d r^{2k+1}\text{d}r = \frac{2}{d^2}
\left[\frac{r^{2k+2}}{2k+2}\right]_0^r = \frac{d^{2k}}{k+1}
\end{multline*}
for $k\in\mathbb N$. Thus 
\begin{multline} \label{D+T}
  \norm{\sum_{n=0}^\infty \lambda^n (D+T)^n}_2^2 = \frac{1}{|\lambda|^2} \left[
\text{exp}\left(\sum_{k=0}^\infty 
|\lambda|^{2(k+1)}\frac{d^{2k}}{k+1}\right)-1\right] \\
= \frac{1}{|\lambda|^2} \text{exp}\left(\frac{1}{d^2}\left(-\log(1-d^2|\lambda|^2)\right)\right) \\
= \frac{1}{|\lambda|^2}\left[(1-d^2|\lambda|^2)^{-\frac{1}{d^2}}-1\right].
\end{multline}
If instead $D+cT$ is a $\text{DT}(\mu_D,c)$ operator with $\mu_D$ being the uniform distribution on a disc of radius $d$ then 
\begin{equation*}
  D+cT = c(D'+T)
\end{equation*}
where $D'$ now has the uniform distribution on $B(0,\frac{d}{c})$, so from (\ref{D+T}) we obtain
\begin{multline}
  \label{D+cT}
   \norm{\sum_{n=0}^\infty \lambda^n (D+cT)^n}_2^2 \\
=\norm{\sum_{n=0}^\infty (c\lambda)^n (D'+T)^n}_2^2 
= \frac{1}{c^2|\lambda|^2}\left[\left(1-d^2|\lambda|^2\right)^{-\frac{c^2}{d^2}}-1\right].
\end{multline}

\begin{lemma} \label{lemma:4.6a} Let $a>b>0$ and let $S=\sqrt{a}T_1+\sqrt{b}T_2^*$ where $T_1$ and $T_2$ are two $\mathcal D$-free quasidiagonal DT-operators. Then
  \begin{equation*}
    \norm{\sum_{n=0}^\infty \lambda^nS^n}_2^2 = 
 \frac{1}{|\lambda|^2}
\frac{\text{e}^{(a-b)|\lambda|^2}-1}{a-b\text{e}^{(a-b)|\lambda|^2}},~~ |\lambda|< \frac{1}{\norm{S}^2}.
  \end{equation*}
\end{lemma}
\begin{proof}
Let $F_n(x) = E_\mathcal D((S^*)^n S^n)$ for $n\in\mathbb N$ and $x\in [0,1]$. For $t< \frac{1}{\norm{S}^2}$ define the $\mathcal D$-valued function  
\begin{equation} \label{F(t,x)}
  F(t,x) = \sum_{n=0}^\infty F_n(x) t^n.
\end{equation}
By Speicher's cumulant formula we have by $\mathcal D$-Gaussianity of $S$ that
\begin{multline*}
  F_n =E_\mathcal D((S^*)^nS^n) = \sum_{\pi\in \text{NC}(2n)} \kappa_\pi^\mathcal D\left( (S^*)^{\otimes_\mathcal B n}\otimes_\mathcal B S^{\otimes_\mathcal B n}\right) \\
= \kappa_2^\mathcal D\left(S^*\otimes_\mathcal B E_\mathcal D((S^*)^{n-1}S^{n-1}) S\right) \\
= (aL^*+bL)(E_\mathcal D((S^*)^{n-1}S^{n-1})) = (aL^*+bL)(F_{n-1}),
\end{multline*}
so we get the following recursive algorithm for determining the $F_n$'s.
\begin{equation*}
  \begin{cases}
    F_0(x) =1  & \\
    F_n(x) =  aL^*(F_{n-1})(x) + bL(F_{n-1})(x), & x\in [0,1] 
  \end{cases},
\end{equation*}
where $L^*(f):x\mapsto \int_0^x f(y)\text{d}y$ and $L(f):x\mapsto \int_x^1 f(y)\text{d}y$. Observe that
\begin{eqnarray*}
  \frac{\text{d}}{\text{d}x}L(f)(x) = -f(x) & \text{ and } &  \frac{\text{d}}{\text{d}x}L^*(f)(x) = f(x),
\end{eqnarray*}
and that 
\begin{equation*}
  F_n(0) = aL^*(F_{n-1})(0) + bL(F_{n-1})(0) = b\int_0^1 F_{n-1}(x)\text{d}x = b\tau(F_{n-1}) 
\end{equation*}
for $n\geq 1$. Using (\ref{F(t,x)}) we have the following differential equation and initial condition in $x$   
\begin{equation*}
  \begin{cases}
    \frac{\text{d}}{\text{d}x} F(t,x) = (a-b)t F(t,x), & x\in [0,1] \\
    F(t,0) = f(t),
  \end{cases}
\end{equation*}
where the function $f$ is given by 
\begin{multline*}
  f(t) = F(t,0) = \sum_{n=0}^\infty F_n(0)t^n \\
 = 1 + \sum_{n=1}^\infty \left( aL^*(F_{n-1})(0) + bL(F_{n-1})(0)\right)t^n \\
= 1 + b\sum_{n=1}^\infty \left(\int_0^1 F_{n-1}(x)\text{d}x\right) t^n \\
= 1 + bt\int_0^1 \left(\sum_{n=1}^\infty F_{n-1}(x)t^{n-1}\right)\text{d}x \\
=1+ bt \tau(F(t,\cdot)) 
\end{multline*}

We thus have the unique solution
\begin{equation} \label{F}
  F(t,x) = f(t)\text{e}^{(a-b)tx},
\end{equation}
where we can now use (\ref{F}) and the initial condition to find the 
function $f$. 
\begin{multline*}
  f(t) = 1 +bt\int_0^1 F(t,x)\text{d}x \\
  = 1 + bt\left[\frac{f(t)}{(a-b)t}\text{e}^{(a-b)tx}\right]_0^1 
 = 1+bf(t)\frac{\left( \text{e}^{(a-b)t} -1\right)}{a-b}.
\end{multline*}
Hence 
\begin{equation*}
  f(t) = \frac{a-b}{a-b\text{e}^{(a-b)t}}
\end{equation*}
so that
\begin{equation*}
  F(t,x) = \frac{(a-b)\text{e}^{(a-b)tx}}{a-b\text{e}^{(a-b)t}}.
\end{equation*}
Now observe that
\begin{multline*}
    \norm{\sum_{n=0}^\infty \lambda^nS^n}_2^2 = \tau\left(F(|\lambda|^2,x)\right) \\
= \int_0^1 F(|\lambda|^2,x)\text{d}x = \frac{1}{|\lambda|^2}
\frac{\text{e}^{(a-b)|\lambda|^2}-1}{a-b\text{e}^{(a-b)|\lambda|^2}}
\end{multline*}
\end{proof}

\begin{theorem}
  The operator $T+\sqrt{\epsilon}Y$ is not a \rm{DT}-operator.
\end{theorem}
\begin{proof}
  By substituting $a=1+\epsilon$ and $b=\epsilon$ in lemma \ref{lemma:4.6a} we have
  \begin{equation}
    \label{eq:4.24a}
    \norm{\sum_{n=0}^\infty \lambda^n(T+\sqrt{\epsilon}Y)^n}^2_2 = \frac{1}{|\lambda|^2}\frac{\text{e}^{|\lambda|^2}-1}{1+\epsilon -\epsilon \text{e}^{|\lambda|^2}}
  \end{equation}
for all $\lambda$ in a neighborhood of $0$.
If $T+\sqrt{\epsilon}Y$ is a DT-operator, then by Theorem \ref{theorem:4.3a} and (\ref{D+cT}), there exists a $c>0$, such that when $d= \log(1+\frac{1}{\epsilon})^{-\frac{1}{2}}$
\begin{equation}
  \label{eq:4.25a}
   \norm{\sum_{n=0}^\infty \lambda^n(T+\sqrt{\epsilon}Y)^n}^2_2= \frac{1}{c^2|\lambda|^2} 
\left((1-d^2|\lambda|^2)^{-\frac{c^2}{d^2}} 
-1\right)
\end{equation}
for all $\lambda$ in a neighborhood of $0$. Consider the two analytic functions,
\begin{eqnarray*}
  f(s)&=& \frac{\text{e}^{s}-1}{1+\epsilon -\epsilon \text{e}^{s}}, \\
 g(s) & =& \frac{1}{c^2} 
\left((1-d^2 s)^{-\frac{c^2}{d^2}} 
-1\right)
\end{eqnarray*}
which are both defined in the complex disc $U=B(0,\log(1+\frac{1}{\epsilon})^{-\frac{1}{2}})$. By (\ref{eq:4.24a}) and (\ref{eq:4.25a}) $f(s)=g(s)$ for s in some real interval of the form $(0,\delta)$ and hence $f(s)=g(s)$ for all $s\in U$. Moreover $f$ has a meromorphic extension to the full complex plane with a simple pole at $s_0=\log(1+\tfrac{1}{\epsilon})$. Hence $g$ also has a meromorphic extension to the full complex plane with a simple pole at $\log(1+\tfrac{1}{\epsilon}) = d^{-2}$. This implies $c=d$. In this case 
\begin{equation*}
  g(s) = \frac{1}{d^2}\left((1-d^2 s)^{-1}-1\right)
\end{equation*}
which is analytic in $\mathbb C\setminus \{s_0\}$. However $f$ has infinitely many poles, namely
\begin{equation*}
  s_p = \log\left(1+\frac{1}{\epsilon}\right)+p2\pi,~~p\in \mathbb Z.
\end{equation*}
Since the meromorphic extensions of $f$ and $g$ must coincide, we have reached a contradiction. Therefore $T+\sqrt{\epsilon}Y$ is not a DT-operator.
\end{proof}

\section{\'Sniady's moment formulas. The case $k=2$.} \label{Sniadydiffsection}

Let $k\in \mathbb N$ be fixed, and let $(P_{k,n})_{n=0}^\infty$ be the sequence of polynomials defined recursively by
  \begin{equation}
    \label{eq:pkn}
    \begin{cases}
      P_{k,n}(x)=1, & \\
      P_{k,n}^{(k)}(x) = P_{k,n-1}(x+1), & n=1,2,\ldots \\
      P_{k,n}(0)=P_{k,n}^{(1)}(0)= \cdots P_{k,n}^{(k-1)}(0)=0, & n=1,2,\ldots
    \end{cases},
  \end{equation}
where $P_{k,n}^{(l)}$ denotes the $l$'th derivative of $P_{k,n}$. As in the previous sections, $T$ denotes the quasinilpotent $DT$ operator. \'Sniady's main results from \cite{Sniady} are:

\begin{theorem}\cite[Theorem 5 and Theorem 7]{Sniady} 
\label{Th5.1}
  \begin{enumerate}[(a)]
  \item For all $k,n\in\mathbb N$:
\begin{equation}
  \label{eq:5.4}
  E_\mathcal D\left(((T^*)^k T^k)^n\right)(x) = P_{k,n}(x), ~~~x\in [0,1].
\end{equation}
\item For all $k,n\in\mathbb N$:
  \begin{equation}
    \label{eq:5.5}
    \tau\left(((T^*)^k T^k)^n\right)= \frac{n^{nk}}{(nk+1)!}
  \end{equation}
  \end{enumerate}  
\end{theorem}

Actually \'Sniady considers $E_\mathcal D((T^k (T^*)^k)^n)$ instead of $E_\mathcal D(((T^*)^k T^k)^n)$, but it is easily seen, that Theorem \ref{Th5.1} (a) is equivalent to \cite[Theorem 5]{Sniady}, by the simple change of variable $x\mapsto 1-x$.

\'Sniady's proof of Theorem \ref{Th5.1} is a very technical combinatorial proof. In this and the following section we will give an analytical proof of Theorem \ref{Th5.1} based on Voiculescu's $\mathcal R$-transform with amalgamation.

As in \cite[(2.11)]{invsub} we put
\begin{equation*}
  \rho(z) = - W_0(-z), ~~~~~~ z\in\mathbb C\setminus [\tfrac{1}{\text{e}},\infty ),
\end{equation*}
where $W_0$ is the principal branch of Lambert's W-function. Then $\rho$ is the principal branch of the inverse function of $z\mapsto z \text{e}^{-z}$. We shall need the following result from \cite[Prop. 4.2]{invsub}. 

\begin{lemma}\cite[Prop. 4.2]{invsub}\label{unique}
  Let $(P_{k,n})_{n=0}^\infty$ be a sequence of polynomials given by (\ref{eq:pkn}). Put for $s\in\mathbb C$, $|s|<\tfrac{1}{\rm{e}}$ and $j=1,\ldots,k$
\begin{eqnarray}
    \alpha_j(s) &=& \rho\left( s\rm{e}^{\rm{i}\frac{2\pi \text{j}}{\text{k}}}\right), \\ 
\gamma_j(s) &=&\begin{cases}
\Pi_{l\neq j} \frac{\alpha_l(s)}{\alpha_l(s)-\alpha_j(s)}, &  0<|s|<\frac{1}{\rm{e}}\\
    \frac{1}{k}, & s=0.
  \end{cases}
\end{eqnarray}
Then 
\begin{equation}
  \label{eq:5.10}
  \sum_{n=0}^\infty (ks)^{nk}P_{k,n}(x) = \sum_{j=1}^k \gamma_j(s)\rm{e}^{\text{k}\alpha_j(\text{s})\text{x}}
\end{equation}
for all $x\in \mathbb R$ and all $s\in B(0,\tfrac{1}{\rm{e}})$.
\end{lemma}

%\begin{proof}
%  By \cite[Prop 4.2]{invsub} the formula (\ref{eq:5.10}) holds for $\phi_n(x) =P_{k,n}(x)$. The uniqueness part of lemma \ref{unique} follows by applying for each fixed $x\in \mathbb R$ the uniqueness of the coefficients in the power series expansion  of the analytic function
%  \begin{equation*}
%    s\mapsto \sum_{j=1}^k \gamma_j(s)\text{e}^{k\alpha_j(s)x}
%  \end{equation*}
%for $s\in \mathbb C\setminus [\tfrac{1}{\text{e}},\infty)$.
%\end{proof}

The case $k=1$ of theorem \ref{Th5.1} is the special case $\lambda=0$ of theorem \ref{Pr3.1}. To illustrate our method of proof of theorem \ref{Th5.1} for $k\geq 2$, we first consider the case $k=2$.

Define $\tilde{T}\in M_4(\mathcal A)$ by
\begin{equation*}
  \tilde{T} = 
  \begin{pmatrix}
    0 & 0 & 0 & T^* \\
    T & 0 & 0 & 0   \\
    0 & T & 0 & 0   \\
    0 & 0 & T^* & 0
  \end{pmatrix}.
\end{equation*}
Then $\lVert\tilde{T}\rVert=\norm{T}=\sqrt{\text{e}}$. (cf. \cite[Corollary 8.11]{DT})
For $\mu\in\mathbb C$, $|\mu|< \tfrac{1}{\text{e}}$ we let $z=z(\mu)$, denote the Cauchy transform of $\tilde{T}$ at $\tilde{\mu} = \mu 1_{M_4(\mathcal A)}$ wrt. amalgamation over $M_4(\mathcal D)$ i.e.
\begin{equation*}
  z =E_\mathcal D\left((\tilde{\mu}-\tilde{T})^{-1}\right).
\end{equation*}
Clearly
\begin{equation}
  \label{eq:5.11}
  (\tilde{\mu}-\tilde{T})^{-1}= \sum_{n=0}^\infty \mu^{-n-1}\tilde{T}^n  =
\left(\sum_{n=0}^3 \mu^{-n-1}\tilde{T}^n\right)\left(\sum_{n=0}^\infty \mu^{-4n} \tilde{T}^{4n}\right).
\end{equation}
By direct computation
\begin{equation*}
  \tilde{T}^2=
  \begin{pmatrix}
    0 & 0 &(T^*)^2 & 0 \\
    0 & 0 & 0 &TT^* \\
    T^2 & 0 & 0 & 0 \\
    0 & T^*T & 0 & 0 
  \end{pmatrix},
\end{equation*}
\begin{equation*}
  ~~\tilde{T}^3 = 
  \begin{pmatrix}
    0 & (T^*)^2T & 0 & 0 \\
    0 & 0 & T(T^*)^2 & 0 \\
    0 & 0 & 0 & T^2T^*   \\
    T^* T^2 &0 &0 &0
  \end{pmatrix}
\end{equation*}
and 
\begin{equation*}
  \tilde{T}^4= 
  \begin{pmatrix}
    (T^*)^2 T^2 & 0 & 0 & 0 \\
    0 & T(T^*)^2 T & 0 & 0  \\
    0 & 0 & T^2(T^*)^2 & 0 \\
    0 & 0 & 0 & T^* T^2 T^* 
  \end{pmatrix}.
\end{equation*}
Hence using the fact that the expectation $E_\mathcal D$ of a monomial in $T$ and $T^*$ vanishes unless $T$ and $T^*$ occur the same number of times, we get from (\ref{eq:5.11}) that $z$ is of the form
\begin{equation}
  \label{eq:5.12}
  z = 
  \begin{pmatrix}
    z_{11} & 0 & 0 & 0 \\
    0 & z_{22} & 0 & z_{24} \\
    0 & 0 & z_{33} & 0 \\
    0 & z_{42} & 0 & z_{44}
  \end{pmatrix}
\end{equation}
where $z_{11},z_{22},z_{24},z_{33},z_{42},z_{44}\in \mathcal D$ are given by
\begin{eqnarray*}
  z_{11} & = & \mu^{-1} E_\mathcal D((1 - \mu^{-4} (T^*)^2 T^2)^{-1}), \\
  z_{22} & = & \mu^{-1} E_\mathcal D((1 - \mu^{-4} T(T^*)^2 T)^{-1}), \\
  z_{33} & = & \mu^{-1} E_\mathcal D((1 - \mu^{-4} T^2 (T^*)^2)^{-1}), \\
  z_{44} & = & \mu^{-1} E_\mathcal D((1 - \mu^{-4} T^* T^2 T^*)^{-1}), \\
  z_{24} & = & \mu^{-3} E_\mathcal D(T(1 - \mu^{-4} (T^*)^2 T^2)^{-1}T^*), \\
  z_{42} & = & \mu^{-3} E_\mathcal D(T^*(1 - \mu^{-4} T^2 (T^*)^2)^{-1}T).
\end{eqnarray*}
For the last 2 identities, we have used, that
\begin{equation*}
  A(1-\eta BA)^{-1} = (1-\eta AB)^{-1} A
\end{equation*}
for $A,B\in \mathcal A$ and $\eta\in\mathbb C$ whenever both sides of this equality are welldefined.

By lemma \ref{RCtheorem}, we know, that there exists a $\delta>0$ such that when $w\in M_4(\mathcal D)_{\text{inv}}$ and $\mu\in\mathbb C$ satisfies $\norm{w}<\delta$, $|\mu|> \tfrac{1}{\delta}$ and
\begin{equation} \label{eq:5.13}
  \mathcal R_{\tilde{T}}^{M_4(\mathcal D)}(w) + w^{-1} = \mu 1_{M_4(\mathcal A)}
\end{equation}
then $w= E_{M_4(\mathcal D)}((\tilde{\mu}-\tilde{T})^{-1}) = z$. In particular 
\begin{equation*}
  w_{11}= z_{11} = \mu^{-1}((1-\mu^{-4}(T^*)^2 T^2)^{-1}),
\end{equation*}
Hence, if we can find a suitable solution to (\ref{eq:5.12}) for all $\mu\in\mathbb C$ in a neighborhood of $\infty$, we can find $E_\mathcal D(((T^*)^2 T^2)^n)$ for $n=1,2,\ldots$ by determining the power series expansion of $w_{11}$ as a function of $\mu^{-1}$.

Since $(T,T^*)$ is a $\mathcal D$-Gaussian pair by \cite[Appendix]{invsub} it follows from lemma \ref{tensorcumulant} that 
\begin{equation*}
  \kappa_n^{M_4(\mathcal D)}((m_1\otimes a_1)\otimes_{M_4(\mathcal D)} \cdots 
\otimes_{M_4(\mathcal D)}(m_n\otimes a_n)) = 0 
\end{equation*}
when $n\neq 2$, $m_1,m_2,\ldots , m_n \in M_4(\mathbb C )$ and $a_1,a_2,\ldots, a_n \in\{T,T^*\}$. By definition
\begin{equation*}
  \tilde{T} =(e_{21}+e_{32})\otimes T + (e_{43}+e_{14})\otimes T^*
\end{equation*}
so by linearity of $\kappa_n^{M_4(\mathcal D)}$, it follows that 
\begin{equation*}
  \kappa_n^{M_4(\mathcal D)}(\tilde{T}\otimes_{M_4(\mathcal D)}\cdots \otimes_{M_4(\mathcal D)}\tilde{T})=0
\end{equation*}
when $n\neq 2$ i.e. $\tilde{T}$ is $M_4(\mathcal D)$-Gaussian. 

Hence using (\ref{R-D-Gaussian}) we get
\begin{multline*}
\mathcal R_{\tilde{T}}^{M_4(\mathcal D)}(w) = \kappa_2^{M_4(\mathcal D)}(\tilde{T}\otimes_{M_4(\mathcal D)} w \tilde{T}) 
= E_{M_4(\mathcal D)}\left(\tilde{T}w\tilde{T}\right) \\ = 
E_{M_4(\mathcal D)}\left( 
  \begin{pmatrix}
  T^*w_{42}T & 0 & T^*w_{44}T^* & 0 \\
  0 & 0 & 0 & Tw_{11}T^* \\
  Tw_{22}T & 0 & Tw_{24}T^* & 0 \\
 0 & T^*w_{33}T & 0 & 0
  \end{pmatrix}
\right) 
\end{multline*}
for $w= (w_{ij})_{i,j=1,\ldots,4}\in M_4(\mathcal D)$.

Since $E_\mathcal D(TfT)=E_\mathcal D(T^*fT^*)=0$, and $E_\mathcal D(T^*fT) = L^*(f)$, $E_\mathcal D(TfT^*)=L(f)$ for $f\in L^\infty([0,1])$, we have:
\begin{equation*}
  \mathcal R_{\tilde{T}}^{M_4(\mathcal D)}(w) = 
  \begin{pmatrix}
    L^*(w_{42}) & 0 & 0 & 0 \\
    0 & 0 & 0 & L(w_{11}) \\
    0 & 0 & L(w_{24}) & 0 \\
    0 & L^*(w_{33}) & 0 & 0
  \end{pmatrix}
\end{equation*}
for $w\in M_4(\mathcal D)$. By (\ref{eq:5.12}) we only have to consider $w$ of the form 
\begin{equation} \label{eq:5.14}
  w = 
  \begin{pmatrix}
    w_{11} & 0 & 0 & 0 \\
    0 & w_{22} & 0 & w_{24} \\
    0 & 0 & w_{33} & 0 \\
    0 & w_{42} & 0 & w_{44}
  \end{pmatrix}.
\end{equation}
For $w\in M_4(\mathcal D)_{\text{inv}}$ of the form (\ref{eq:5.14}), (\ref{eq:5.13}) reduces to the three equations
\begin{equation} \label{eq:5.15}
  \begin{cases}
    L^*(w_{42}) + \frac{1}{w_{11}} = \mu 1_\mathcal D & \\
      \begin{pmatrix}
        0 & L(w_{11}) \\
        L^*(w_{33}) & 0
      \end{pmatrix}
+ 
\begin{pmatrix}
  w_{22} & w_{24} \\
  w_{42} & w_{44}
\end{pmatrix}^{-1} = \mu 1_{M_2(\mathcal D)} & \\
L(w_{24}) + \frac{1}{w_{33}} = \mu 1_\mathcal D &
  \end{cases}.
\end{equation}

\begin{definition}
  Let $f\in C([0,1])$. We call $(f^{(-n)})_{n=1}^l$ for the succesive antiderivatives of $f$ if
  \begin{equation*}
    \frac{\text{d}}{\text{d}x}(f^{(-n)})= f^{(1-n)} \text{ for } n=2,3,\ldots, l
  \end{equation*}
and
\begin{equation*}
  \frac{\text{d}}{\text{d}x}(f^{(-1)})=f.
\end{equation*}
\end{definition}

\begin{lemma} \label{lemma:5.4}
Let $f\in C^2([0,1])$ and let $f^{(-1)}$ and $f^{(-2)}$ be the succesive antiderivatives of $f$ for which
  \begin{enumerate}
\item[(i)]   \label{eq:i} $f^{(-1)}(1)=0, ~~~~f^{(-2)}(1)=\mu^3$.
 \end{enumerate}
Assume further, that
\begin{enumerate}
\item[(ii)] $f(0)=\mu^{-1}$ and $f^{(1)}(0)=0$. \label{eq:ii} \\
\item[(iii)] \label{eq:iii} For all $x\in [0,1]$,
  \begin{eqnarray*}
    f(x)& \neq&  0 \\
    \left|\begin{matrix}
      f^{(-1)}(x) & f(x) \\
      f(x) & f^{(1)}(x)
    \end{matrix}\right| &\neq &0
  \end{eqnarray*}
while 
\begin{equation*}
  \left|
    \begin{matrix}
      f^{(-2)}(x) & f^{(-1)}(x) & f(x) \\
      f^{(-1)}(x) & f(x)  & f^{(1)}(x) \\
      f(x)& f^{(1)}(x) & f^{(2)}(x) 
    \end{matrix}\right|=0
\end{equation*}.
\end{enumerate}
Then $w_{11}, w_{22}, w_{33}, w_{44},w_{24}, w_{42}\in C([0,1])$ given by
\begin{equation}
  \label{eq:5.16}
  \begin{cases}
    w_{11} = f \\
w_{22}=w_{44} = -\frac{1}{\mu}\frac{
  \left|\begin{matrix}
     f^{(-1)} & f \\
     f & f^{(1)}
  \end{matrix}\right|
}{f^2} \\
~ \\
w_{24} = \frac{1}{\mu^2} \frac{ f^{(-1)}\left|\begin{matrix}
     f^{(-1)} & f \\
     f & f^{(1)}
  \end{matrix}\right|}{f^2} \\
~ \\
 w_{42} = \frac{f^{(1)}}{f^2} \\
~ \\
w_{33} = \mu^2\frac{ f\left|\begin{matrix}
     f & f^{(1)} \\
     f^{(1)} & f^{(2)}
  \end{matrix}\right|}{ \left|\begin{matrix}
     f^{(-1)} & f \\
     f & f^{(1)}
  \end{matrix}\right|^2}
  \end{cases}
\end{equation}
is a solution to (\ref{eq:5.15}). Moreover
\begin{equation} \label{eq:5.17}
   \left|\begin{matrix}
     w_{22} & w_{24} \\
     w_{42} & w_{44}
  \end{matrix}\right| = -\frac{1}{\mu^2}\frac{ \left|\begin{matrix}
     f^{(-1)} & f \\
     f & f^{(1)}
  \end{matrix}\right|}{f^2}
\end{equation}
and
\begin{equation}
  \label{eq:5.18}
  \begin{cases}
    L(w_{11}) = -f^{(-1)} \\
~ \\
    L(w_{24}) = \mu - \frac{1}{\mu^2}\frac{ \left|\begin{matrix}
     f^{(-2)} & f^{(-1)} \\
     f^{(-1)} & f
  \end{matrix}\right|}{f} \\ 
~ \\
    L^*(w_{42})= \mu-\frac{1}{f} \\
~ \\
    L^*(w_{33}) = -\mu^2 \frac{f^{(1)}}{ \left|\begin{matrix}
     f^{(-1)} & f \\
     f & f^{(1)}
  \end{matrix}\right|}
  \end{cases}.
\end{equation}
\end{lemma}
\begin{proof}
  Assume $w_{11}, w_{22}, w_{33}, w_{44}, w_{24}, w_{42}$ is given by (\ref{eq:5.16}). Then (\ref{eq:5.17}) follows immediately. Note that for $f\in C([0,1])$, the functions $g=L(f)$ and $h=L^*(f)$ are characterized by
  \begin{eqnarray*}
    g^{(1)} = -f & \text{ and } & g(1)=0 \\
    h^{(1)} = f  & \text{ and } & h(0)=0.
  \end{eqnarray*}
Hence (\ref{eq:5.18}) is equivalent to (\ref{eq:5.19}) and (\ref{eq:5.20}) 
below.
\begin{equation} \label{eq:5.19}
  \begin{cases}
    \frac{\text{d}}{\text{d}x}f^{(-1)} = w_{11} \\
~ \\
    \frac{\text{d}}{\text{d}x}\Bigl(\frac{1}{\mu^2}\frac{\left|\begin{matrix}
     f^{(-2)} & f^{(-1)} \\
     f^{(-1)} & f
  \end{matrix}\right|}{f}\Bigr) = w_{24} \\
~ \\
    \frac{\text{d}}{\text{d}x}\left(-\frac{1}{f}\right) = w_{42} \\
~ \\
     \frac{\text{d}}{\text{d}x}\Bigl(-\mu^2 \frac{f^{(1)}}{\left|\begin{matrix}
     f^{(-1)} & f \\
     f & f^{(1)}
  \end{matrix}\right|} \Bigr) = w_{33}
  \end{cases}
\end{equation}
\begin{equation}
  \label{eq:5.20}
  \begin{cases}
    f^{(-1)}(1)=0, & \frac{\left|\begin{matrix}
     f^{(-2)}(1) & f^{(-1)}(1) \\
     f^{(-1)}(1) & f(1)
  \end{matrix}\right|}{f(1)}= \mu^3 \\
~ \\
\frac{1}{f(0)}= \mu, & f^{(1)}(0)=0
  \end{cases}.
\end{equation}
Now, (\ref{eq:5.20}) is trivial from (i) and (ii). Next we prove (\ref{eq:5.19}): Clearly
\begin{eqnarray*}
  \frac{\text{d}}{\text{d}x}f^{(-1)}=f=w_{11} & \text{ and } &
 \frac{\text{d}}{\text{d}x}(-\frac{1}{f}) =\frac{f^{(1)}}{f^2} = w_{42}.
\end{eqnarray*}
Moreover
\begin{multline}
  \frac{\text{d}}{\text{d}x}\Bigl(\frac{\left|\begin{matrix}
     f^{(-2)} & f^{(-1)} \\
     f^{(-1)} & f
  \end{matrix}\right|}{f}\Bigr) \\ = 
\frac{ f\left|\begin{matrix}
     f^{(-2)} & f \\
     f^{(-1)} & f^{(1)}
  \end{matrix}\right| - f^{(1)}\left|\begin{matrix}
     f^{(-2)} & f^{(-1)} \\
     f^{(-1)} & f
  \end{matrix}\right|}{f^2} 
= \frac{f^{(-1)}\left|\begin{matrix}
     f^{(-1)} & f \\
     f & f^{(1)}
  \end{matrix}\right|}{f^2}= \mu^2 w_{24}
\end{multline}
and
\begin{multline*}
   \frac{\text{d}}{\text{d}x}\Bigl(\frac{f^{(1)}}{\left|\begin{matrix}
     f^{(-1)} & f \\
     f & f^{(1)}
  \end{matrix}\right|}\Bigr) = \frac{\left|\begin{matrix}
     f^{(-1)} & f \\
     f & f^{(1)}
  \end{matrix}\right|f^{(2)} - \left|\begin{matrix}
     f^{(-1)} & f \\
     f^{(1)} & f^{(2)}
  \end{matrix}\right|f^{(1)}}{\left|\begin{matrix}
     f^{(-1)} & f \\
     f & f^{(1)}
  \end{matrix}\right|^2} \\
= -\frac{f\left|\begin{matrix}
     f & f^{(1)} \\
     f^{(1)} & f^{(2)}
  \end{matrix}\right| }{\left|\begin{matrix}
     f^{(-1)} & f \\
     f & f^{(1)}
  \end{matrix}\right|^2}= -\frac{1}{\mu^2} w_{33}.
\end{multline*}
Hence (\ref{eq:5.19}) holds. It remains to be proved that $w_{11}, w_{22}, w_{33}, w_{44}, w_{24},w_{42}$ is a solution to (\ref{eq:5.15}). By (\ref{eq:5.16}) and (\ref{eq:5.18}), we have
\begin{equation*}
  L^*(w_{42}) + \frac{1}{w_{11}} = \left(\mu - \frac{1}{f}\right) + \frac{1}{f} = \mu.
\end{equation*}
Moreover by (\ref{eq:5.16}) and (\ref{eq:5.17})
\begin{multline*}
  \begin{pmatrix}
    w_{22} & w_{24} \\
    w_{42} & w_{44}
  \end{pmatrix}^{-1} = \frac{1}{w_{22}w_{44}-w_{24}w_{42}}
  \begin{pmatrix}
    w_{44} & -w_{24} \\
    -w_{42} & w_{22}
  \end{pmatrix}  \\ 
= 
\begin{pmatrix}
\mu & f^{(-1)} \\
 & \\
\mu^2\frac{f^{(1)}}{\left|
    \begin{matrix}
    f^{(-1)} & f \\
    f & f^{(1)}
    \end{matrix}
\right|}  & \mu
\end{pmatrix}
\end{multline*}
which proves that the first and the second inequality in (\ref{eq:5.15}).

By (\ref{eq:5.16}) and  (\ref{eq:5.18}),
\begin{equation*}
  w_{33}(\mu-L(w_{24})) = \frac{\left|
    \begin{matrix}
    f & f^{(1)} \\
    f^{(1)} & f^{(2)}
    \end{matrix}
\right|\left|
    \begin{matrix}
    f^{(-2)} & f^{(-1)} \\
    f^{(-1)} & f
    \end{matrix}
\right|}{\left|
    \begin{matrix}
    f^{(-1)} & f \\
    f & f^{(1)}
    \end{matrix}
\right|^2}  = 1 + \frac{\sigma}{\left|
    \begin{matrix}
    f^{(-1)} & f \\
    f & f^{(1)}
    \end{matrix}
\right|^2}
\end{equation*}
where 
\begin{equation*}
  \sigma = \left|
    \begin{matrix}
    f & f^{(1)} \\
    f^{(1)} & f^{(2)}
    \end{matrix}
\right|\left|
    \begin{matrix}
    f^{(-2)} & f^{(-1)} \\
    f^{(-1)} & f
    \end{matrix}
\right| - \left|
    \begin{matrix}
    f^{(-1)} & f \\
    f & f^{(1)}
    \end{matrix}
\right|^2 = f \left|
    \begin{matrix}
    f^{(-2)} & f^{(-1)} & f  \\
    f^{(-1)} &  f       & f^{(1)} \\
    f & f^{(1)} & f^{(2)}
    \end{matrix}
\right|.
\end{equation*}
Hence by (iii), $\sigma=0$. Therefore $w_{33}(x) \neq 0$ for all $x\in [0,1]$ and $w_{33}^{-1} = \mu - L(w_{24})$, proving the last equality in (\ref{eq:5.15}).
\end{proof}

\begin{lemma} \label{lemma:5.5}
  Let $\alpha_j(s), \gamma_j(s)$ for $j=1,2$ be as in lemma \ref{unique} for $k=2$, i.e. $\alpha_1(0)=\alpha_2(0)=0$, $\gamma_1(0)=\gamma_2(0)=\tfrac{1}{2}$ and for $0<|s|< \rm{e}^{-1}$:
  \begin{eqnarray*}
    \alpha_1(s)=\rho(s), & \text{       } &  \alpha_2(s)=\rho(-s), \\
   \gamma_1(s) = \frac{\alpha_1(s)}{\alpha_1(s) -\alpha_2(s)}, & \text{      } &
\gamma_2(s) = \frac{\alpha_2(s)}{\alpha_2(s) -\alpha_1(s)}.
  \end{eqnarray*}
Let $\mu\in\mathbb C$, $|\mu|>\sqrt{\text{e}}$, put $s=\tfrac{1}{2}\mu^{-2}$ and
\begin{eqnarray}
  \label{eq:5.21} f(x) & = & \frac{1}{\mu}\left(\sum_{j=1}^2 \gamma_j(s)\text{e}^{2\alpha_j(s)x}\right), \text{  } x\in\mathbb R \\
  \label{eq:5.22} f^{(-1)}(x) & = & \frac{1}{2\mu}\left(\sum_{j=1}^2 
\frac{\gamma_j(s)}{\alpha_j(s)}\text{e}^{2\alpha_j(s)x}\right), \text{  } x\in\mathbb R \\
  \label{eq:5.23} f^{(-2)}(x) & = & \frac{1}{4\mu}\left(\sum_{j=1}^2 
\frac{\gamma_j(s)}{\alpha_j(s)^2}\text{e}^{2\alpha_j(s)x}\right), \text{  } x\in\mathbb R 
\end{eqnarray}
Then
\begin{enumerate}[(i)]
\item $f^{(-1)}, f^{(-2)}$ are succesively antiderivatives of $f$,
  \begin{equation}
    \label{eq:5.24}
    f^{(-1)}(1)=0, \text{  } f^{(-2)}(1) = \mu^3
  \end{equation}
and
  \begin{equation}
    \label{eq:5.25}
    f(0)=\mu^{-1}, \text{  } f^{(1)}(0) = 0.
  \end{equation}
\item The following asymptotic formulas holds for $|\mu|\to \infty$:
  \begin{eqnarray*}
    f^{(-2)}(x) &=& \mu^3 + \mathcal O(\mu^{-1}) \\
 f^{(-1)}(x) &=& (x-1)\mu^{-1} + \mathcal O(\mu^{-5}) \\
 f(x) &=& \mu^{-1} + \mathcal O(\mu^{-5}) \\
 f^{(1)}(x) &=& x\mu^{-5} + \mathcal O(\mu^{-9}) \\
 f^{(2)}(x) &=& x\mu^{-5} + \mathcal O(\mu^{-9}) 
  \end{eqnarray*}
where the error estimates holds uniformly in $x$ on a compact subset in $\mathbb R$.
\item There exists $\mu_0\geq \sqrt{\rm{e}}$ such that the restriction of $f$ to $[0,1]$ satisfies all the conditions in lemma \ref{lemma:5.4}, when $|\mu| > \mu_0$.
\end{enumerate}
\end{lemma}

\begin{proof}
Clearly $f^{(-1)}$ and $f^{(-2)}$ are succesively antiderivatives of $f$ and
      \begin{eqnarray*}
        f(0) & =& \frac{1}{\mu}\sum_{j=1}^2 \gamma_j(s) = \frac{1}{\mu} \\
        f^{(1)}(0) & = & \frac{2}{\mu}\sum_{j=1}^2 \alpha_j(s)\gamma_j(s) = 0.
      \end{eqnarray*}
To prove (\ref{eq:5.24}), note first, that since $\rho:\mathbb C\setminus [\frac{1}{\text{e}},\infty)\to \mathbb C$ is a branch of the inverse function of 
$z\mapsto z\text{e}^{-z}$, we have
\begin{equation*}
  \rho(w)\text{e}^{-\rho(w)}=w, \text{   } |w|<\frac{1}{\text{e}}
\end{equation*}
and therefore
\begin{equation*}
  \text{e}^{2\alpha_j(s)}= \frac{\alpha_j(s)^2}{s^2}, \text{  }j=1,2.
\end{equation*}
Since $s^2=\frac{1}{4}\mu^{-4}$, it follows that
\begin{eqnarray}
  \label{eq:5.26} f^{(-2)}(x+1) & = & \mu^4f(x), \text{   } x\in\mathbb R \\
 \label{eq:5.27} f^{(-1)}(x+1) & = & \mu^4f^{(1)}(x), \text{   } x\in\mathbb R \\
 \label{eq:5.28} f(x+1) & = & \mu^4f^{(2)}(x), \text{   } x\in\mathbb R.
\end{eqnarray}
In particular
\begin{eqnarray*}
  f^{(-2)}(1) &=& \mu^4f(0) = \mu^3 \\
  f^{(-1)}(1) &=& \mu^4f^{(1)}(0) = 0.
\end{eqnarray*}
By the proof of \cite[Prop. 4.2]{invsub}, $\alpha_j(s)$ and $\rho_j(s)$ are continuous functions of $s\in B(0,\tfrac{1}{\text{e}})$. Hence, regarding $f$ as a function of $\mu$, 
  \begin{equation*}
    \lim_{|\mu|\to\infty} (\mu f(x)) = \sum_{j=1}^2 \gamma_j(0)\text{e}^{2\alpha_j(0)x}=1
  \end{equation*}
where the limit holds uniformly in $x$ on compact subsets of $\mathbb R$. Hence by (\ref{eq:5.28}) $f^{(2)}(x)=\mathcal O(\mu^{-5})$ as $|\mu| \to \infty$ uniformly in $x$ on compact subsets of $\mathbb R$. By (\ref{eq:5.25}),
\begin{eqnarray}
  \label{eq:5.29} f^{(1)}(x) & = & \int_0^x f^{(2)}(t)\text{d}t \\
  \label{eq:5.30} f(x) & = & \mu^{-1} + \int_0^x f^{(1)}(t)\text{d}t
\end{eqnarray}
which implies, that $f^{(1)}(x)= \mathcal O(\mu^{-5})$ and
\begin{equation}
  f(x) = \mu^{-1} + \mathcal O(\mu^{-5}) \label{eq:5.31}
\end{equation}
uniformly in $x$ on compact subsets of $\mathbb R$.

Using again (\ref{eq:5.28}), (\ref{eq:5.29}) and (\ref{eq:5.30}), we get
\begin{eqnarray*}
  f^{(2)}(x) & = & \mu^{-5} + \mathcal O(\mu^{-9}) \\
  f^{(1)}(x) & = & x\mu^{-5} + \mathcal O(\mu^{-9}).
\end{eqnarray*}
By (\ref{eq:5.24})
\begin{eqnarray*}
  f^{(-1)}(x) &= & \int_1^x f(t)\text{d}t \\
  f^{(-2)}(x) & = & \mu^3 +\int_1^x f^{(-1)}(t) \text{d}t.
\end{eqnarray*}
Hence by (\ref{eq:5.31}),
\begin{eqnarray*}
  f^{(-1)}(x) &= & (x-1)\mu^{-1} + \mathcal O(\mu^{-5}) \\
  f^{(-2)}(x) &= & \mu^{3} + \mathcal O(\mu^{-1})
\end{eqnarray*}
where all estimates holds uniformly on compact subsets of $\mathbb R$.
 This proves (ii). 
\item By (i), $f^{(-1)}, f^{(-2)}$ coinside with the succesive antiderivatives of $f$ considered in lemma \ref{lemma:5.4} and $f(0)=\mu^{-1}, f^{(1)}(0) = 0$.

Moreover, by (ii),
\begin{eqnarray*}
  f(x) & = & \mu^{-1} + \mathcal O(\mu^{-5}) \\
 \left|
   \begin{matrix}
   f^{(-1)}(x) & f(x) \\
   f(x) & f^{(1)}(x)
   \end{matrix}
\right| &=&  \mu^{-2} + \mathcal O(\mu^{-6})
\end{eqnarray*}
where the error terms holds uniformly in $x\in [0,1]$. Hence there exists $\mu_0\geq \sqrt{\rm{e}}$, such that
\begin{equation*}
  f(x) \neq 0 \text{ and } \left|
   \begin{matrix}
   f^{(-1)}(x) & f(x) \\
   f(x) & f^{(1)}(x)
   \end{matrix}
\right| \neq 0
\end{equation*}
for all $x\in [0,1]$. Moreover by the matrix factorization
\begin{multline}\label{eq:5.29a}
  \begin{pmatrix}
    f^{(-2)}(x) & f^{(-1)}(x) & f(x) \\
  f^{(-1)}(x) & f(x) & f^{(1)}(x) \\
  f(x) & f^{(1)}(x) & f^{(2)}(x) \\
  \end{pmatrix} \\ =  
  \begin{spmatrix}
    1 & 1 \\
2\alpha_1(s) & 2\alpha_2(s) \\
4\alpha_1(s)^2 & 4\alpha_2(s)^2
  \end{spmatrix} 
  \begin{spmatrix}
    \frac{\gamma_1(s)}{4\alpha_1(s)^2}\text{e}^{2\alpha_1(s)x} & 0 \\
0 &  \frac{\gamma_2(s)}{4\alpha_2(s)^2}\text{e}^{2\alpha_2(s)x}
  \end{spmatrix}
  \begin{spmatrix}
    1 & 2\alpha_1(s) & 4\alpha_1(s)^2 \\
  1 & 2\alpha_2(s) & 4\alpha_2(s)^2 
  \end{spmatrix}
\end{multline}
it follows, that the matrix on the left hand side has rank less than or equal to 2, i.e.
\begin{equation*}
\left|\begin{matrix}
     f^{(-2)}(x) & f^{(-1)}(x) & f(x) \\
  f^{(-1)}(x) & f(x) & f^{(1)}(x) \\
  f(x) & f^{(1)}(x) & f^{(2)}(x) \\
\end{matrix}\right| = 0
\end{equation*}
for $x\in [0,1]$. Hence $f$ satisfies all the conditions in lemma \ref{lemma:5.4}, when $|\mu|>\mu_0$. 
\end{proof}

\begin{proof}[Proof of Theorem \ref{Th5.1} in the case $k = 2$:]
  
By lemma \ref{RCtheorem} there exists a $\delta > 0$, such that when $w\in M_4(\mathcal D)_{\text{inv}}$ and $\mu\in\mathbb C$ satisfies $\norm{w}<\delta, |\mu|>\frac{1}{\delta}$ and
\begin{equation}
  \label{eq:5.32}
  \mathcal R_{\tilde{T}}^{M_4(\mathcal D)}(w) + w^{-1} = \mu 1_{M_4(\mathcal D)}
\end{equation}
then $w= E_\mathcal D((\tilde{\mu} - \tilde{T})^{-1})$. In particular
\begin{equation}
  \label{eq:5.33}
  w_{11} = \mu^{-1} E_\mathcal D((1 - \mu^{-4}(T^*)^2 T^2)^{-1}).
\end{equation}
Let $\mu\in\mathbb C, |\mu| > \sqrt{\text{e}}$, put $s= \frac{1}{2}\mu^{-2}$ and
\begin{equation*}
  f(x)  =  \frac{1}{\mu} \left( \sum_{j=1}^2 \gamma_j(s) \text{e}^{2\alpha_j(s)x}\right) 
\end{equation*}
for $x\in [0,1]$ as in lemma \ref{lemma:5.5}. By lemma \ref{lemma:5.5} (iii) there exists a $\mu_0> \sqrt{\text{e}}$, such that when $|\mu|> \mu_0$, then $f$ satisfies all the requirements af lemma \ref{lemma:5.4}. Hence by lemma \ref{lemma:5.4}, the matrix $w\in M_4(\mathcal D)$ given by (\ref{eq:5.14}) and (\ref{eq:5.16}) is a solution to (\ref{eq:5.32}). Moreover by the asymptotic formulas in lemma \ref{lemma:5.5} (ii),
\begin{eqnarray*}
   \left|
   \begin{matrix}
   f^{(-2)}(x) & f^{(-1)}(x) \\
   f^{(-1)}(x) & f(x)
   \end{matrix}
\right| &=& \mu^2 + \mathcal O(\mu^{-2}), \\
  \left|
   \begin{matrix}
   f^{(-1)}(x) & f(x) \\
   f(x) & f'(x)
   \end{matrix}
\right|&=& -\mu^{-2} + \mathcal O(\mu^{-6}), \\
  \left|
   \begin{matrix}
   f(x) & f'(x) \\
   f'(x) & f''(x)
   \end{matrix}
\right|&=& \mu^{-6} + \mathcal O(\mu^{-10}).
\end{eqnarray*}
Hence by (\ref{eq:5.16}) and the asymptotic formulas for $f^{(-1)}, f$ and $f'$, we have
\begin{eqnarray*}
  w_{11} & = & \mu^{-1} + \mathcal O(\mu^{-5}), \\
w_{22} & = & w_{44} = \mu^{-1} + \mathcal O(\mu^{-5}), \\
w_{24} & = & (1-x)\mu^{-3} + \mathcal O(\mu^{-3}), \\
w_{42} & = & x\mu^{-3} + \mathcal O(\mu^{-3}), \\
w_{33} & = & \mu^{-1} + \mathcal O(\mu^{-5}), \\
\end{eqnarray*}
where all the error estimates holds uniformly in $x\in [0,1]$. Hence, there exists $\mu_1 \geq \max\{\mu_0, \tfrac{1}{\delta}\}$, such that when $|\mu|> \mu_1$ then $\norm{w}< \delta$, and hence 
\begin{equation*}
  w = E_{M_4(\mathcal D)} ((\tilde{\mu}-\tilde{T})^{-1}).
\end{equation*}
By (\ref{eq:5.16}), $w_{11} = f$. Hence by (\ref{eq:5.33}) and (\ref{eq:5.21})
\begin{equation*}
  E_\mathcal D((1-\mu^{-4}(T^*)^2 T^2)^{-1})(x) = \mu f(x) = \sum_{j=1}^2 \gamma_j(s) \text{e}^{2\alpha_j(s)x}
\end{equation*}
where $s=\frac{1}{2}\mu^{-2}$, i.e. for $|s|<\frac{1}{2}\mu_1^{-2}$, 
\begin{equation*}
  E_\mathcal D((1-(2s)^2(T^*)^2 T^2)^{-1})(x) = \sum_{j=1}^2 \gamma_j(s) \text{e}^{2\alpha_j(s)x}
\end{equation*}
and therefore 
\begin{equation} \label{eq:5.31n}
  \sum_{j=0}^\infty (2s)^{2n} E_\mathcal D(((T^*)^2 T^2)^n)(x) = \sum_{j=1}^2 \gamma_j(s) \text{e}^{2\alpha_j(s)x}.
\end{equation}
Hence by lemma \ref{unique} and by the uniqueness of the power series expansions of analytic functions, we have 
\begin{equation*}
  E_\mathcal D (((T^*)^2 T^2)^n)(x) = P_{2,n}(x)
\end{equation*}
for $n\in\mathbb N$ and $x\in [0,1]$. This proves theorem \ref{Th5.1}(a) in the case $k=2$. Theorem 5.1 (b) also follows from (\ref{eq:5.31n}) by integrating the right hand side of (\ref{eq:5.31n}) from 0 to 1 with respect to $x$ (cf. \cite[remark 4.3]{invsub}).
\end{proof} 

\section{\'Sniady's moment formulas. The general case.} \label{Sniadygeneralcase}

The above proof of Theorem \ref{Th5.1} in the case $k=2$ can fairly easily be generalized to all $k\geq 2$ (Recall that the case $k=1$ is contained in theorem \ref{Pr3.1}).

Let $k\geq 2$ and define  $\tilde{T}\in M_{2k}(\mathcal A)$ by
  \begin{equation*}
    \tilde{T} =\sum_{j=1}^k (T\otimes e_{j+1,j} + T^*\otimes e_{k+j+1,k+j})
  \end{equation*}
where the indices are computed modulo $2k$, such that $e_{2k+1,2k}= e_{1,2k}$. For $\mu\in\mathbb C, |\mu| < \frac{1}{\sqrt{\text{e}}}$, we put $\tilde{\mu} = \mu 1_{2k}$ and 
\begin{equation*}
  z= z(\mu) = E_{M_{2k}(\mathcal D)}((\tilde{\mu}-\tilde{T})^{-1}).
\end{equation*}
Then only the diagonal entries $z_{11}, \ldots, z_{2k,2k}$ and the off-diagonal entries 
$z_{2,2k}, z_{3,2k-1}, \ldots, z_{2k,2}$ can be non-zero. Moreover,
\begin{equation*}
  z_{11} = \mu^{-1} E_\mathcal D((1-\mu^{-2k} (T^*)^k T^k)^{-1}).
\end{equation*}
The operator $\tilde{T}$ is $M_{2k}(\mathcal D)$-Gaussian, and repeating the arguments for $k=2$, we get that for $w\in M_{2k}(\mathcal D)$, the matrix 
\begin{equation}\label{eq:6.1}
u=\mathcal R_{\tilde{T}}^{M_{2k}(\mathcal D)}(w)
\end{equation}
can have at most $2k$ non-zero entries, namely the entries
\begin{eqnarray} 
\nonumber  u_{11} & = & L^*(w_{2k,2}) \\
\nonumber  u_{2k,2} & = & L^*(w_{2k-1,3}) \\
\nonumber  \vdots ~~~ & & ~~~~\vdots \\
\label{eq:6.2}  u_{k+2,k} & = & L^*(w_{k+1,k+1}) \\
\nonumber u_{k+1,k+1} & = & L(w_{k,k+2}) \\
\nonumber u_{k,k+2} & = & L(w_{k-1,k+3}) \\
\nonumber \vdots~~~ & & ~~~~\vdots \\
\nonumber u_{2,2k} & = & L(w_{1,1}).
\end{eqnarray}
By lemma \ref{RCtheorem} there exists a $\delta>0$ (depending on $k$), such that if $w\in M_{2k}(\mathcal D)_{\text{inv}}, \norm{w}<\delta, \mu\in\mathbb C, |\mu| > \frac{1}{\delta}$ and
  \begin{equation} \label{eq:6.3}
    \mathcal R_{\tilde{T}}^{M_{2k}(\mathcal D)}(w) + w^{-1} = \mu 1_{M_{2k}(\mathcal D)},
  \end{equation}
then
\begin{equation*}
  w=z=E_{M_{2k}(\mathcal D)}((\tilde{\mu}-\tilde{T})^{-1}).
\end{equation*}
In particular
\begin{equation*}
  w_{11}= \mu^{-1}E_\mathcal D((1- \mu^{-2k}(T^*)^{k}T^{k})^{-1}).
\end{equation*}
Next we construct an explicit solution to (\ref{eq:6.3}). By the above remarks on $z$, it is sufficient to consider those $w\in M_{2k}(\mathcal D)_{\text{inv}}$ for which only the entries $z_{11}, \ldots, z_{2k,2k}$ and $z_{2,2k}, z_{3,2k-1}, \ldots, z_{2k,2}$ can be non-zero. For such $w$, (\ref{eq:6.3}) can by 
(\ref{eq:6.1}) and (\ref{eq:6.2}) be reduced to the $k+1$ identities:
\begin{equation}  \label{eq:6.4} 
\begin{cases}
L^*(w_{2k,2})+ \frac{1}{w_{11}} = \mu 1_\mathcal D \\
~\\
\begin{spmatrix}
  0 & L(w_{j+1,2k+1-j}) \\
 L^*(w_{2k-1-j,j+3}) & 0 
\end{spmatrix} + 
\begin{spmatrix}
  w_{2+j,2+j} & w_{2+j,2k-j} \\
  w_{2k-j,2+j} & w_{2k-j,2k-j}
\end{spmatrix}^{-1} = \mu 1_{M_2(\mathcal D)}, \\
 ~~~~~~~~~~~~~~~~~~~~~~~~~~~~~~~~~~~~~~~~~~~~~
~~~~~~~~~~~~~~~~~~~~~j=0,1,\ldots, k-2, \\
L(w_{k,k+2}) + \frac{1}{w_{k+1,k+1}} = \mu 1_\mathcal D.
\end{cases}
\end{equation}

\begin{definition} \label{definition:6.1}
  For $j\in \mathbb N\cup\{0\}$ and $g\in C^{2j+2}$, we let $\Delta_j(g)$ denote the determinant
  \begin{equation}
    \label{eq:6.5}
   \Delta_j(g)= \left|\begin{smallmatrix}
g       & g^{(1)}  & \adots & g^{(j)} \\
g^{(1)} & \adots   & \adots & \adots  \\
\adots  & \adots   & \adots & g^{(2j-1)} \\
g^{(j)} & \adots   & g^{(2j-1)}& g^{(2j)}
\end{smallmatrix}\right|.
  \end{equation}
In particular $\Delta_0(g) =g$.
\end{definition}

\begin{lemma}\label{lemma:5.7}
  Let $g\in C^{2j+2}(\mathbb R)$ and $j\in \mathbb N$. Then
  \begin{equation}
    \label{eq:6.6}
    \Delta_j(g^{(2)})\Delta_j(g) - \Delta_j(g^{(1)})^2 = \Delta_{j-1}(g^{(2)}) \Delta_{j+1}(g)
  \end{equation}
\end{lemma}
and
  \begin{equation}
    \label{eq:6.7}
 \Delta_{j-1}(g^{(2)}) \frac{\text{d}}{\text{d}x}\left(
 \Delta_{j}(g)\right)-  \Delta_{j}(g)\frac{\text{d}}{\text{d}x} 
\left(
\Delta_{j-1}(g^{(2)})\right) =  \Delta_{j-1}(g^{(1)}) 
\Delta_{j}(g^{(1)}).
  \end{equation}

The proof of lemma \ref{lemma:5.7} relies on elementary matrix manipulations and is contained in lemma \ref{Hankel} of appendix \ref{app1}. More specifically (\ref{eq:6.6}) is a direct consequence of (\ref{Hankela}) from lemma \ref{Hankel}, and (\ref{eq:6.7}) follows from (\ref{Hankelb}) of lemma \ref{Hankel} by using the elementary fact that:
\begin{equation*}
   \frac{\text{d}}{\text{d}x}\left(\Delta_j(g)\right)= \left|\begin{smallmatrix}
g       & g^{(1)}  & \adots & g^{(j)} \\
g^{(1)} & \adots   & \adots & \adots  \\
\adots  & \adots   & \adots & g^{(2j-1)} \\
g^{(j-1)} & \adots   & g^{(2j-2)}& g^{(2j-1)} \\
g^{(j+1)} & \adots   & g^{(2j)}& g^{(2j+1)} 
\end{smallmatrix}\right|,
\end{equation*}
that is, differentiating (\ref{eq:6.5}) is the same as differentiating the last row of (\ref{eq:6.5}).

The next two lemmas are the generalizations of lemma \ref{lemma:5.4} and lemma \ref{lemma:5.5} to arbitrary $k\geq 2$. 
\begin{lemma} \label{lemma:6.3}
  Let $f\in C^k([0,1])$ and let $(f^{(-j)})^k_{j=1}$ be the antiderivatives of $f$ for which,
  \begin{enumerate}[(i)]
  \item  \label{6.3(i)} \begin{equation*}
    f^{(-j)}(1) =
    \begin{cases}
      0, & 1 \leq j \leq k-1, \\
      \mu^{2k-1}, & j=k.
    \end{cases}
 \end{equation*}
\item \label{6.3(ii)} Assume further that 
\begin{equation*}
  f(0) = \mu^{-1} \text{ and } f^{(-j)}(0)=0 \text{ for } 1\leq j\leq 
k-1.
\end{equation*}
\item \label{6.3(iii)} For all $x\in [0,1]$,
\begin{equation*}
  \Delta_j(f^{(-j)})(x) \neq 0, \text{ for } j=0\ldots, k-1
\end{equation*}
and
\begin{equation*}
  \Delta_k(f^{(-k)})(x) =0
\end{equation*}
\end{enumerate}
Then the set of $4k-2$ functions listed in (\ref{eq:6.8}), (\ref{eq:6.9}) and (\ref{eq:6.10}) below is a solution to (\ref{eq:6.4}).
\begin{equation}
  \label{eq:6.8}
  \begin{cases}
    w_{11} = f & \\
~\\
    w_{22}=w_{2k,2k} = -\frac{1}{\mu}\frac{\Delta_1(f^{(-1)})}{f^2} \\
~\\
    w_{2,2k} = \frac{1}{\mu^2} \frac{f^{(-1)}\Delta_1(f^{(-1)})}{f^2}  \\
~\\
w_{2k,2} = \frac{f^{(1)}}{f^2}
  \end{cases}.
\end{equation}
For $j=1,\ldots,k-2$
\begin{equation}
  \label{eq:6.9}
  \begin{cases}
    w_{j+2,j+2} = w_{2k-j,2k-j} = 
-\frac{1}{\mu}\frac{\Delta_{j-1}(f^{(1-j)})\Delta_{j+1}(f^{(-1-j)})}{\Delta_j(f^{(-j)})^2}  \\
~\\
w_{j+2,2k-j}= \frac{1}{\mu^{2j+2}} \frac{\Delta_j(f^{(-1-j)})\Delta_{j+1}(f^{(-1-j)})}{\Delta_j(f^{(-j)})^2}  \\
~\\
w_{2k-j,j+2} = \mu^{2j} \frac{\Delta_{j-1}(f^{(1-j)})\Delta_{j}(f^{(1-j)})}{\Delta_j(f^{(-j)})^2}
  \end{cases}.
\end{equation}
\begin{equation}
  \label{eq:6.10}
  w_{k+1,k+1} = \mu^{2k+2} \frac{\Delta_{k-2}(f^{(2-k)})\Delta_{k-1}(f^{(2-k)})}{\Delta_{k-1}(f^{(1-k)})^2} 
\end{equation}
Moreover for $j=0,\ldots, k-2$
\begin{equation}
  \label{eq:6.11}
  \left|
    \begin{matrix}
    w_{j+2,j+2} & w_{j+2,2k-2} \\
    w_{2k-j,j+2}& w_{2k-j,2k-j}
    \end{matrix}
\right|= \frac{1}{\mu} w_{j+2,j+2}
\end{equation}
and
\begin{equation}
  \label{eq:6.12}
  \begin{cases}
    L(w_{11}) = -f^{(-1)} & \\
~\\
    L(w_{j+2,2k-j}) = -\frac{1}{\mu^{2j+2}}\frac{\Delta_{j+1}(f^{(-2-j)})}{\Delta_j(f^{(-j)})}, &  0\leq j\leq k-3 \\
~\\
    L(w_{k,k+2}) = \mu - \frac{1}{\mu^{2k-2}}\frac{\Delta_{k-1}(f^{(-k)})}{\Delta_{k-2}(f^{(2-k)})}
  \end{cases}
\end{equation}
\begin{equation}
  \label{eq:6.13}
  \begin{cases}
    L^*(w_{2k,2}) = \mu- \frac{1}{f} & \\
~\\
    L^*(w_{2k-j,2+j}) = -\mu^{2j}\frac{\Delta_{j-1}(f^{(-2-j)})}{\Delta_j(f^{(-j)})}, & 1\leq j \leq k-2 \\ 
~\\
L^*(w_{k+1,k+1}) = -\mu^{2k-2} \frac{\Delta_{k-2}(f^{(3-k)})}{\Delta_{k-1}(f^{(1-k)})}
  \end{cases}.
\end{equation}
\end{lemma}

\begin{proof}
  Let $w_{11},w_{22},\ldots, w_{kk}, w_{2,2k},w_{3,2k-1},\ldots,w_{2k,2}$ be given by (\ref{eq:6.8}), (\ref{eq:6.9}) and (\ref{eq:6.10}). Then for $1\leq j \leq k-2$ the left hand side of (\ref{eq:6.11}) is equal to 
  \begin{equation*}
    -\frac{1}{\mu^2}\frac{\Delta_{j-1}(f^{(1-j)}) 
\Delta_{j+1}(f^{(-1-j)})A}{\Delta_j(f^{(-j)})^4},
  \end{equation*}
where $A= \Delta_{j-1}(f^{(1-j)})\Delta_{j+1}(f^{(-1-j)})-
\Delta_{j}(f^{(1-j)})\Delta_{j}(f^{(-1-j)})$.
 
By applying (\ref{eq:6.6}) to $g=f^{(-1-j)}$ it follows that $A=-\Delta_j(f^{(-j)})^2$, which proves (\ref{eq:6.11}) for $1\leq j \leq k-2$. The case $j=0$ of (\ref{eq:6.11}) follows immediately from (\ref{eq:6.8}).

The proofs of (\ref{eq:6.12}) and \ref{eq:6.13}) can be obtained exactly as in the case $k=2$ provided the following two identities holds:
For $j=0,\ldots,k-2$:
\begin{equation}
  \label{eq:6.14}
\frac{\text{d}}{\text{d}x} \left(\frac{\Delta_{j+1}(f^{(-2-j)})}{\Delta_j(f^{(-j)})}\right)
=
\frac{\Delta_{j}(f^{(-1-j)})\Delta_{j+1}(f^{(-1-j)})}{\Delta_j(f^{(-j)})^2}
\end{equation}
For $j=1,\ldots,k-1$:
\begin{equation}
  \label{eq:6.15}
  \frac{\text{d}}{\text{d}x} \left(\frac{\Delta_{j-1}(f^{(2-j)})}{\Delta_j(f^{(-j)})}\right)
=
\frac{\Delta_{j-1}(f^{(1-j)})\Delta_{j}(f^{(1-j)})}{\Delta_j(f^{(-j)})^2}
\end{equation}
However (\ref{eq:6.14}) follows from (\ref{eq:6.7}) with $g=f^{(-2-j)}$ after changing $j$ in (\ref{eq:6.7}) to $j+1$. In the same way (\ref{eq:6.15}) follows from (\ref{eq:6.7}) with $g= f^{(-j)}$ and $j$ unchanged.
It remaims to be proved, that $w_{11}, \ldots, w_{kk},w_{2,2k},\ldots,w_{2k,2}$ form a solution to (\ref{eq:6.4}). The proof of the first 2 identities in (\ref{eq:6.4}) is exactly the same as in the case $k=2$. Let us check the next $k-2$ identities in (\ref{eq:6.4}) i.e.
\begin{multline}
  \label{eq:6.16}
  \begin{pmatrix}
    0 & L(w_{j+1,2k+1-j}) \\
 L^*(w_{2k-1-j,j+3})& 0 
  \end{pmatrix} \\ + 
  \begin{pmatrix}
    w_{2+j,2+j} & w_{2+j,2k-j} \\
    w_{2k-j,2+j} & w_{2k-j,2k-j}
  \end{pmatrix}^{-1} = \mu 1_{M_2(\mathcal D)}
\end{multline}
for $j=1,\ldots, k-2$. By (\ref{eq:6.11}) and the fact that $w_{2+j,2+j}= w_{2k-j,2k-j}$ (cf. (\ref{eq:6.8})) we have 
\begin{equation*}
    \begin{pmatrix}
    w_{2+j,2+j} & w_{2+j,2k-j} \\
    w_{2k-j,2+j} & w_{2k-j,2k-j}
  \end{pmatrix}^{-1}= 
  \begin{pmatrix}
    \mu 1_\mathcal D & \beta \\
    \gamma & \mu 1_\mathcal D
  \end{pmatrix},
\end{equation*}
where
\begin{equation*}
  \beta = -\mu \frac{w_{2+j,2k-j}}{w_{2+j,2+j}}= \frac{1}{\mu^{2j}}
\frac{\Delta_j(f^{(-1-j)})}{\Delta_{j-1}(f^{(1-j)})}
\end{equation*}
and 
\begin{equation*}
  \gamma = -\mu \frac{w_{2k-j,2+j}}{w_{2+j,2+j}}= \mu^{2j+2}
\frac{\Delta_j(f^{(1-j)})}{\Delta_{j+1}(f^{(-1-j)})}.
\end{equation*}
Hence by (\ref{eq:6.12}) and (\ref{eq:6.13})
\begin{equation*}
  \beta = -L(w_{j+1,2k-j+1}) \text{ and } \gamma = -L^*(w_{2k-1-j,j+3})
\end{equation*}
for $j=1,\ldots, k-2$. This proves (\ref{eq:6.16}). Observe next that by (\ref{eq:6.10}) and (\ref{eq:6.12})
\begin{multline*}
  w_{k+1,k+1}(\mu-L(w_{k,k+2})) = \frac{\Delta_{k-1}(f^{(2-k)})\Delta_{k-1}(f^{(-k)})}{\Delta_{k-1}(f^{(1-k)})^2} \\ =
1+ \frac{\sigma}{\Delta_{k-1}(f^{(1-k)})^2},
\end{multline*}
where 
\begin{equation*}
  \sigma = \Delta_{k-1}(f^{(2-k)})\Delta_{k-1}(f^{(-k)})- \Delta_{k-1}(f^{(1-k)})^2.
\end{equation*}

By (\ref{eq:6.6}) and the assumptions (iii) in lemma \ref{lemma:6.3}
\begin{equation*}
  \sigma = \Delta_{k-2}(f^{(2-k)})\Delta_{k}(f^{(-k)})=0.
\end{equation*}
Hence $w_{k+1,k+1}(\mu-L(w_{k,k+2}))=1$, which proves the last equality in (\ref{eq:6.4}). This completes the proof of lemma \ref{lemma:6.3}.
\end{proof}

\begin{lemma} \label{lemma:6.4}
  Let $k\in \mathbb N, k\geq 2$ and let $\alpha_j(s), \gamma_j(s)$ for $j=1,\ldots,k$ and $0< |s| < \tfrac{1}{\rm{e}}$ be as in lemma \ref{unique}. Let $\mu\in \mathbb C$, $|\mu|>\sqrt{\rm{e}}$, put $s=\frac{1}{k}\mu^{-2}$ and
  \begin{equation}
    \label{eq:6.17}
    \begin{cases}
      f(x) = \frac{1}{\mu}\left(\sum_{\nu=1}^k \gamma_\nu(s) 
\rm{e}^{\text{k}\alpha_\nu(\text{s})\text{x}} \right), & x\in \mathbb R \\
 f^{(-j)}(x) = \frac{1}{\mu k^j}\left(\sum_{\nu=1}^k \frac{\gamma_\nu(s)}{\alpha_\nu(s)^j} 
\rm{e}^{\text{k}\alpha_\nu(\text{s})\text{x}} \right), & x\in \mathbb R, j=1,\ldots, k 
    \end{cases}.
  \end{equation}
Then 
\begin{enumerate}[(i)]
\item $(f^{(-j)})_{j=1}^k$ are succesive antiderivatives of $f$. Moreover
  \begin{equation}
    \label{eq:6.18}
    \begin{cases}
      f^{(-j)}(1)=0, & 1\leq j \leq k-1 \\
      f^{(-k)}(1) = \mu^{2k-1}
    \end{cases}
  \end{equation}
and 
\begin{equation}
  \label{eq:6.19}
  \begin{cases}
    f(0)=\mu^{-1} \\
    f^{(j)}(0) =0, & 1 \leq j \leq k-1
  \end{cases}.
\end{equation}
\item The following asymptotic formulas holds for $|\mu|\to \infty$
  \begin{equation}
    \label{eq:6.20}
    \begin{cases}
      f^{(-k)}(x) = \mu^{2k-1} + \mathcal O(\mu^{-1}) \\
      f^{(-j)}(x) = \frac{1}{j!}(x-1)^j \mu^{-1} + \mathcal O(\mu^{-2k-1}), & 1\leq j\leq k-1 \\
 f(x) = \mu^{-1} + \mathcal O(\mu^{-2k-1}) \\
 f^{(j)}(x) = \frac{1}{j!}x^j \mu^{-2k-1} + \mathcal O(\mu^{-4k-1}), & 1\leq j\leq k-1 \\
f^{(k)}(x) = \mu^{-2k-1} +\mathcal O(\mu^{-4k-1})
    \end{cases},
  \end{equation}
where the error estimates holds uniformly in $x$ on compact subsets of $\mathbb R$.
\item There exists a $\mu_0 \geq \sqrt{\rm{e}}$, such that the restriction of $f$ to $[0,1]$ satisfies all the conditions in lemma \ref{lemma:6.3}, when $|\mu| > \mu_0$.
\end{enumerate}
\end{lemma}

\begin{proof}
  From the proof of \cite[Prop. 4.2]{invsub}, we know that $\alpha_j(s)$ and $\gamma_j(s)$ are analytic functions of $s\in B(0,\tfrac{1}{\text{e}})$. Moreover by \cite[Prop. 4.1]{DT}
  \begin{equation}
    \label{eq:6.21}
    \begin{cases}
      \sum_{\nu=1}^k \gamma_\nu(s) =1 \\
      \sum_{\nu=1}^k \gamma_\nu(s)\alpha_\mu(s)^j =1, & j=1,\ldots, k-1
    \end{cases}.
  \end{equation}
Moreover, since $\alpha_j(s) = \rho\bigl(\text{e}^{\text{i}\tfrac{2\pi j}{k}s}\bigr)$, where $\rho$ satisfies
\begin{equation*}
  \rho(w)\text{e}^{-\rho(w)}= w  \text{ for } |w|< \frac{1}{\text{e}}
\end{equation*}
we have $\bigl(\alpha_\nu(s)\text{e}^{-\alpha_\nu(s)}\bigr)^k= s^k $ and therefore
\begin{equation} \label{eq:6.22}
  \text{e}^{k\alpha_\nu(s)} = \frac{s^k}{(\alpha_\nu(s))^k}
\end{equation}
for $\nu=1,\ldots,k$. Having (\ref{eq:6.21}) and (\ref{eq:6.22}) in mind, the proof of (i) and (ii) in lemma \ref{lemma:6.4} is now a routine generalization of the proof of lemma \ref{lemma:5.5}.
Concerning (iii) in lemma \ref{lemma:6.4}, we have
\begin{equation}
  \label{eq:6.23}
  \begin{cases}
    \Delta_j(f^{(-j)})= \sigma(j)\mu^{-j-1} + \mathcal O(\mu^{-2k-j-1}), & 0,\ldots, k-1, \\
\text{where } \sigma(j)=1 \text{ for } j=0,3 \text{ (mod } 4) \\
\text{and } \sigma(j)=-1 \text{ for } j=1,2 \text{ (mod } 4) 
  \end{cases}
\end{equation}
because the leading term in the determinant $\Delta_j(f^{(-j)})$ comes from the antidiagonal, i.e.
\begin{equation*}
  \Delta_j(f^{(-j)}) = \left|
    \begin{matrix}
      0 & \dots & 0 & f \\
      \vdots & \adots & \adots & 0 \\
      0 & \adots & \adots & \vdots \\
      f & 0 & \dots & 0
    \end{matrix}\right| + \mathcal O(\mu^{-2k-j-1}) = \sigma(j)f^{j+1} + \mathcal O(\mu^{-2k-j-1})
\end{equation*}
since the matrix in question has size $j+1$. 
Hence $\Delta_j(f^{(-j)})(x)\neq 0$ for $x\in [0,1]$ and $0\leq j\leq k-1$, when $|\mu |$ is sufficiently large. Moreover $\Delta_k(f^{(-k)})=0$ for $x\in [0,1]$, because in analogy with (\ref{eq:5.29a}), $\Delta_k(f^{(-k)}(x))$ is the determinant of the $(k+1)\times (k+1)$ matrix
\begin{equation*}
  F=(f^{(i+j-k)})_{i,j=0,\ldots,k}
\end{equation*}
which has the factorization $F=ADA^t$, where $A$ is the $(k+1)\times k$ matrix with entries 
\begin{equation*}
 a_{il}=(k\alpha_l(s))^i, ~~i=0,\ldots, k,~~l=1,\ldots,k 
\end{equation*}
and $D$ is the $k\times k$ diagonal matrix, with diagonal entries
\begin{equation*}
  d_{ll} = \frac{\gamma_l(s)}{(k\alpha_l(s))^k}\rm{e}^{\emph{k}\alpha_\emph{l}(\emph{s})},~~ \emph{l}=1,\ldots,\emph{k}. 
\end{equation*}

\end{proof}

\begin{proof}[Proof of Theorem \ref{Th5.1} in the general case]
  Let $\mu_0$ be as in lemma \ref{lemma:6.4}, let $\mu\in\mathbb C, |\mu|>\mu_0$ and put $s=\frac{1}{k}\mu^{-2}$. Put as before
  \begin{equation*}
    f(x) = \frac{1}{\mu}\left(\sum_{\nu=1}^k \gamma_j(s)\text{e}^{k\alpha_j(s)x}\right)
  \end{equation*}
for $x\in [0,1]$, and define $w_{11},w_{22},\ldots,w_{k,k}, w_{2,2k},w_{3,2k-1},\ldots,w_{2k,2}$ by (\ref{eq:6.8}), (\ref{eq:6.9}) and (\ref{eq:6.10}), and put all other entries of $w\in M_{2k}(\mathcal D)$ equal to $0$. Then by lemma \ref{lemma:6.4}, (\ref{eq:6.4}) holds, and therefore
\begin{equation*}
  \mathcal R_{\tilde{T}}^{M_{2k}(\mathcal D)}(w)+ w^{-1} = \mu 1_{M_{2k}(\mathcal D)}.
\end{equation*}
Let $\delta>0$ be chosen according to lemma \ref{RCtheorem}. If we can find a $\mu_1\geq \max\{\mu_0, \frac{1}{\delta}\}$, such that
\begin{equation}
  \label{eq:6.24}
  |\mu|\geq \mu_1 \Rightarrow \norm{w}<\delta
\end{equation}
then $w=E_{M_{2k}(\mathcal D)}((\tilde{\mu}-\tilde{T})^{-1})$. In particular
\begin{equation} \label{eq:6.25n}
  f = w_{11} = \mu^{-1}E_\mathcal D((1- \mu^{-2k}(T^*)^k T^k)^{-1}),
\end{equation}
and the proof of theorem \ref{Th5.1} for $k\geq 2$ can be completed exactly as in the case $k=2$.
By (\ref{eq:6.23})
\begin{equation}
  \label{eq:6.25}
  \begin{cases}
    \Delta_j(f^{(-j)}) = \mathcal O(\mu^{-j-1}), & 0\leq j \leq k-1 \\
     \frac{1}{\Delta_j(f^{(-j)})} = \mathcal O(\mu^{j+1}), & 0\leq j \leq k-1
  \end{cases}
\end{equation}
uniformly in $x\in [0,1]$ for $|\mu|\to \infty$. We claim that
\begin{equation}
  \label{eq:6.26}
  \begin{cases}
    \Delta_j(f^{(-j-1)}) = \mathcal O(\mu^{-j-1}), & 0\leq j \leq k-2\\     \Delta_{k-1}(f^{(-k)}) = \mathcal O(\mu^{k}) \\
 \Delta_j(f^{(1-j)}) = \mathcal O(\mu^{-j-2k-1}), & 0\leq j \leq k-2\\
 \Delta_{k-1}(f^{(2-k)}) = \mathcal O(\mu^{-3k})
  \end{cases}.
\end{equation}
Recall by definition \ref{definition:6.1} that
\begin{equation*}
  \Delta_j(g) = \det\left((g^{(k+l)})_{k,l=0,\ldots,j}\right).
\end{equation*}
Hence for $0\leq j \leq k-2$, $\Delta_j(f^{(-j-1)})$ is the determinant of a $(j+1)\times (j+1)$ matrix, where each entry is equal to one of the functions $f^{(-j-1)}, f^{(-j)}, \ldots, f^{(j-1)}$. 
By (\ref{eq:6.20}) all these functions are of order $\mathcal O(\mu^{-1})$ as $|\mu|\to \infty$. Hence $\Delta_j(f^{(-j-1)})=\mathcal O(\mu^{-j-1})$ proving the first estimate in (\ref{eq:6.26}).
By the same argument, $\Delta_{k-1}(f^{(-k)})$ is the determinant of a $k\times k$ matrix for which the upper left entry is of the order $\mathcal O(\mu^{2k-1})$ and all the other entries are of order $\mathcal O(\mu^{-1})$. Hence $\Delta_{k-1}(f^{(-k)})=\mathcal O(\mu^{2k-1}(\mu^{-1})^{k-1})=\mathcal O(\mu^k)$. Let $0\leq j\leq k-1$. 
Then $\Delta_j(f^{(1-j)})$ is by (\ref{eq:6.20}) a determinant of a $(j+1)\times (j+1)$ matrix $M=(m_{k,l})_{k,l=0,\ldots, j}$ for which
\begin{equation*}
  \begin{cases}
    m_{k,l}=\mathcal O(\mu^{-1}) & \text{ when } k+l < 0 \\
    m_{k,l}=\mathcal O(\mu^{-2k-1}) & \text{ when } k+l \geq 0
  \end{cases}.
\end{equation*}
Hence for any permutation $\pi$ of $\{0,1,\ldots, k\}$ the product 
\begin{equation*}
  m_{0\pi(0)}m_{1\pi(1)} \cdots m_{j\pi(j)}
\end{equation*}
contains at least one factor of order $\mathcal O(\mu^{-2k-1})$. Therefore 
\begin{equation*}
  \Delta_j(f^{(1-j)}) =\det (M) = 
\sum_{\pi\in S_{j+1}}(-1)^{\text{sign}(\pi)} m_{0\pi(0)}m_{1\pi(1)} \cdots m_{k\pi(k)}
\end{equation*}
is of order $\mathcal O(\mu^{-2k-1}(\mu^{-1})^j)=\mathcal O(\mu^{-2k-j-1})$. This proves the last two estimates in (\ref{eq:6.26}). Clearly all estimates holds uniformly in $x\in [0,1]$. Combining (\ref{eq:6.8}), (\ref{eq:6.9}), (\ref{eq:6.10}) and (\ref{eq:6.26}), we get
\begin{equation*}
  \begin{cases}
    w_{l,l} = \mathcal O(\mu^{-1}), & 1\leq l \leq 2k \\
    w_{j+2,2k-j} = \mathcal O(\mu^{-2j-3}), & 0\leq j\leq k-2 \\
    w_{2k-j,j+2} = \mathcal O(\mu^{2j+1-2k}), & 0 \leq j\leq k-2
  \end{cases}.
\end{equation*}
In particular all the entries of w are of size $\mathcal O(\mu^{-1})$ as $|\mu|\to \infty$ uniformly in $x\in [0,1]$. Hence there exists $\mu_1
\geq \max\{\mu_0,\frac{1}{\delta}\}$ such that (\ref{eq:6.24}) holds. Hence by (\ref{eq:6.25n}) we have for $|s|<\frac{1}{k}\mu_1^{-2}$, 
\begin{equation*}
  \sum_{k=0}^\infty (ks)^{nk} E_\mathcal D(((T^*)^k T^k)^n)(x) = 
\sum_{\nu=1}^\infty \gamma_j(s) \rm{e}^{\emph{k}\alpha_j(\emph{s}) \emph{x}}, ~~~\emph{x}\in [0,1].
\end{equation*}
Now Theorem \ref{Th5.1} follows from lemma \ref{unique} and \cite[remark 4.3]{invsub} as in the case $k=2$.
\end{proof}

\appendix

\section{Determinant-identities on Hankel-matrices} \label{app1}

We need the following lemma on Hankel-determinants.
\begin{lemma} \label{Hankel}
  Let $a_{-(n-1)}, a_{-(n-2)}, \ldots, a_{n-1},a_{n}\in \mathbb C$ for
  some $n\in \mathbb N$. Then 
  \begin{enumerate}[(a)]
  \item  \label{Hankela}    \begin{multline*}
    \left| \begin{smallmatrix}
        a_{-(n-3)}    &    a_{-(n-4)} & \adots    &  a_0     \\
        a_{-(n-4)}    & \adots & \adots    & \adots   \\
        \adots & \adots & \adots    & a_{n-4} \\
        a_0    & \adots & a_{n-4} & a_{n-3}
      \end{smallmatrix} \right |
 \left|\begin{smallmatrix}
 a_{-(n-1)}    &    a_{-(n-2)}   &    a_{-(n-3)}    &     a_{-(n-4)} &
 \adots    &  a_0     \\ 
 a_{-(n-2)}    &    a_{-(n-3)}   &    a_{-(n-4)}    &  \adots & \adots
 & \adots   \\ 
 a_{-(n-3)}    &    a_{-(n-4)}   &  \adots   &  \adots & \adots    & a_{n-4} \\ 
 a_{-(n-4)}    &   \adots &  \adots   &  \adots & a_{n-4}  & a_{n-3} \\
 \adots &   \adots &  \adots   & a_{n-4}& a_{n-3}  & a_{n-2} \\
 a_0    &   \adots & a_{n-4}  & a_{n-3}& a_{n-2}  & a_{n-1}
      \end{smallmatrix}\right| \\
=   \left|\begin{smallmatrix}
        a_{-(n-1)}    &    a_{-(n-2)} & \adots    &  a_{-1} \\
        a_{-(n-2)}    & \adots & \adots    & \adots   \\
        \adots & \adots & \adots    & a_{n-4} \\
       a_{-1} & \adots & a_{n-4} & a_{n-3}
      \end{smallmatrix}\right| 
    \left|\begin{smallmatrix}
        a_{-(n-3)}    &    a_{-(n-4)} & \adots    &  a_{1} \\
        a_{-(n-4)}    & \adots & \adots    & \adots   \\
        \adots & \adots & \adots    & a_{n-2} \\
       a_{1} & \adots & a_{n-2} & a_{n-1}
      \end{smallmatrix}\right| 
-  \left|\begin{smallmatrix}
        a_{-(n-2)}    &    a_{-(n-3)} & \adots    &  a_0     \\
        a_{-(n-3)}    & \adots & \adots    & \adots   \\
        \adots & \adots & \adots    & a_{n-3} \\
        a_0    & \adots & a_{n-3} & a_{n-2}
      \end{smallmatrix}\right|^2.
    \end{multline*}
\item \label{Hankelb}
  \begin{multline*}
\left|\begin{smallmatrix}
        a_{-(n-2)}    &    a_{-(n-3)} & \adots    &  a_{1}     \\
        a_{-(n-3)}    & \adots & \adots    & \adots   \\
        \adots & \adots & \adots    & a_{n-1} \\
        a_{1}& \adots & a_{n-1} & a_{n}
      \end{smallmatrix}\right| 
\left|\begin{smallmatrix}
        a_{-(n-2)}    &    a_{-(n-3)} & \adots    &  a_0     \\
        a_{-(n-3)}    & \adots & \adots    & \adots   \\
        \adots & \adots & \adots    & a_{n-3} \\
        a_0    & \adots & a_{n-3}  & a_{n-2}
      \end{smallmatrix}\right| \\
= \left|\begin{smallmatrix}
        a_{-(n-1)}    &    a_{-(n-2)} & \adots    &  a_0     \\
        a_{-(n-2)}    & \adots & \adots    & \adots   \\
        \adots & \adots & \adots    & a_{n-3} \\
       a_{-1} & \cdots & a_{n-3}  & a_{n-2}  \\
       a_{1} & a_{2}& \cdots    & a_{n}
      \end{smallmatrix}\right|
\left|\begin{smallmatrix}
        a_{-(n-3)}    &    a_{-(n-4)} & \adots    &  a_{1}     \\
        a_{-(n-4)}    & \adots & \adots    & \adots   \\
        \adots & \adots & \adots    & a_{n-2} \\
        a_{1}& \adots & a_{n-2} & a_{n-1}
      \end{smallmatrix}\right| \\
- \left|\begin{smallmatrix}
        a_{-(n-1)}    &    a_{-(n-2)} & \adots    &  a_{0}     \\
        a_{-(n-2)}    & \adots & \adots    & \adots   \\
        \adots & \adots & \adots    & a_{n-2} \\
        a_{0}& \adots & a_{n-2} & a_{n-1}
      \end{smallmatrix}\right| 
\left|\begin{smallmatrix}
        a_{-(n-3)}    &    a_{-(n-4)} & \adots    &  a_{1} \\
        a_{-(n-4)}    & \adots & \adots    & \adots   \\
        \adots & \adots & \adots    & a_{n-3} \\
       a_{0} & \cdots & a_{n-3}  & a_{n-2}  \\
       a_{2} & a_{3}& \cdots    & a_{n}
      \end{smallmatrix}\right|.
 \end{multline*}
  \end{enumerate}
\end{lemma}

\begin{proof}
 To prove (\ref{Hankela}) we actually prove the more general equation 
 \begin{multline} \label{matrix1}
 \left|
   \begin{smallmatrix}
    a_{22} & a_{23}    & \cdots & a_{2,n-1} \\
    a_{32} & a_{33}    & \cdots & a_{3,n-1} \\
    \vdots & \vdots    &        & \vdots    \\
 a_{n-1,2} & a_{n-1,3} & \cdots & a_{n-1,n-1}
   \end{smallmatrix}
\right| \left|
   \begin{smallmatrix}
    a_{11} & a_{12}   & a_{13}  & \cdots & a_{1,n}   \\
    a_{21} & a_{22}   & a_{23}  & \cdots & a_{2,n}   \\
    a_{31} & a_{32}   & a_{33}  & \cdots & a_{3,n}   \\
    \vdots & \vdots   & \vdots  &        & \vdots    \\
   a_{n,1} & a_{n,2}  & a_{n,3} & \cdots & a_{n,n}
   \end{smallmatrix}
\right| \\ =
  \left|
   \begin{smallmatrix}
    a_{11} & a_{12}    & \cdots & a_{1,n-1} \\
    a_{21} & a_{22}    & \cdots & a_{2,n-1} \\
    \vdots & \vdots    &        & \vdots    \\
 a_{n-1,1} & a_{n-1,2} & \cdots & a_{n-1,n-1}
   \end{smallmatrix}
\right| \left|
   \begin{smallmatrix}
    a_{22} & a_{23}    & \cdots & a_{2,n} \\
    a_{32} & a_{33}    & \cdots & a_{3,n} \\
    \vdots & \vdots    &        & \vdots    \\
 a_{n,2} & a_{n,3} & \cdots & a_{n,n}
   \end{smallmatrix}
\right| \\-  \left|
   \begin{smallmatrix}
    a_{12} & a_{13}    & \cdots & a_{1,n} \\
    a_{22} & a_{23}    & \cdots & a_{2,n} \\
    \vdots & \vdots    &        & \vdots    \\
 a_{n-1,2} & a_{n-1,3} & \cdots & a_{n-1,n}
   \end{smallmatrix}
\right| \left|
   \begin{smallmatrix}
    a_{21} & a_{22}    & \cdots & a_{2,n-1} \\
    a_{31} & a_{32}    & \cdots & a_{3,n-1} \\
    \vdots & \vdots    &        & \vdots    \\
 a_{n,1} & a_{n,2} & \cdots & a_{n,n-1}
   \end{smallmatrix}
\right|
 \end{multline}
for $a_{ij}\in \mathbb C$ and $i,j\in \{1,\ldots,n\}$.

We first add some zero terms to the left-hand side (LHS) of
(\ref{matrix1}).

\begin{multline*}
\text{LHS}=  \left|
   \begin{smallmatrix}
    a_{22} & \cdots & a_{2,n-1} \\
    \vdots & &  \vdots    \\
 a_{n-1,2} & \cdots & a_{n-1,n-1}
   \end{smallmatrix}
\right| \left|
   \begin{smallmatrix}
    a_{11} & \cdots & a_{1,n}   \\
    \vdots  &        & \vdots    \\
   a_{n,1} & \cdots & a_{n,n}
   \end{smallmatrix}
\right| \\ + \sum_{k=2}^{n-1}   
\left|\begin{smallmatrix}
    a_{21} & \cdots & a_{2,k-1} & a_{2,k+1} & \cdots & a_{2,n-1} \\
    a_{31} & \cdots & a_{3,k-1} & a_{3,k+1} & \cdots & a_{3,n-1} \\
    \vdots &        &  \vdots   & \vdots    &        & \vdots    \\
  a_{n-1,1} & \cdots & a_{n-1,k-1} & a_{n-1,k+1} & \cdots & a_{n-1,n-1} \\
\end{smallmatrix}\right| 
\left|\begin{smallmatrix}
    a_{12} & \cdots & a_{1,k} & a_{1,k} & a_{1,k+1} & \cdots & a_{2,n-1} \\
    a_{22} & \cdots & a_{2,k} & a_{2,k} & a_{2,k+1} & \cdots & a_{3,n-1} \\
    \vdots &        &  \vdots & \vdots  & \vdots    &        & \vdots    \\
  a_{n,2}  & \cdots & a_{n,k} & a_{n,k} & a_{n,k+1} & \cdots & a_{n-1,n-1} \\
\end{smallmatrix}\right| 
\end{multline*}
We note that the last matrix in the sum is zero because coloumn $k-1$
and $k$ are equal. Now we expand LHS after the $k$'th coloumn of the
second matrix in the $k$'th addent. We get
\begin{multline*}
\text{LHS} = \sum_{j=1}^n (-1)^{1+j}a_{j,1} 
\left|
   \begin{smallmatrix}
    a_{22} & \cdots & a_{2,n-1}   \\
    \vdots  &        & \vdots    \\
   a_{n-1,2} & \cdots & a_{n-1,n-1}
   \end{smallmatrix}
\right| 
\left|
   \begin{smallmatrix}
    a_{12}   & \cdots & a_{1,n}     \\
    \vdots   &        & \vdots     \\
   a_{j-1,2} & \cdots & a_{j-1,n} \\
   a_{j+1,2} & \cdots & a_{j+1,n} \\
     \vdots  &        & \vdots     \\
   a_{n,2}   & \cdots & a_{n,n} \\
  \end{smallmatrix}
\right|
 \\+ \sum_{k=2}^{n-1}\sum_{j=1}^n (-1)^{k+j}a_{j,k}
\left|
   \begin{smallmatrix}
    a_{21}      & \cdots & a_{2,k-1}   & a_{2,k+1}   & \cdots & a_{2,n-1}   \\
    \vdots      &        & \vdots      & \vdots      &        &   \vdots    \\
    a_{n-1,1}   & \cdots & a_{n-1,k-1} & a_{n-1,k+1} & \cdots & a_{n-1,n-1} \\
   \end{smallmatrix}
\right| 
\left|
   \begin{smallmatrix}
    a_{12}   & \cdots & a_{1,n}     \\
    \vdots   &        & \vdots     \\
   a_{j-1,2} & \cdots & a_{j-1,n} \\
   a_{j+1,2} & \cdots & a_{j+1,n} \\
     \vdots  &        & \vdots     \\
   a_{n,2}   & \cdots & a_{n,n} \\
  \end{smallmatrix}
\right|
\end{multline*}
where $j=1$ and $j=n$ means leave out row $1$ and
$n$ respectively. Switching the indices we have
\begin{multline} \label{matrix2}
\text{LHS} = \sum_{j=1}^n \left|
   \begin{smallmatrix}
    a_{12}   & \cdots & a_{1,n}     \\
    \vdots   &        & \vdots     \\
   a_{j-1,2} & \cdots & a_{j-1,n} \\
   a_{j+1,2} & \cdots & a_{j+1,n} \\
     \vdots  &        & \vdots     \\
   a_{n,2}   & \cdots & a_{n,n} \\
  \end{smallmatrix}
\right| \left( (-1)^{1+j}a_{j,1}\left|
   \begin{smallmatrix}
    a_{22} & \cdots & a_{2,n-1}   \\
    \vdots  &        & \vdots    \\
   a_{n-1,2} & \cdots & a_{n-1,n-1}
   \end{smallmatrix}
\right|  \right. \\
\left. \sum_{k=2}^{n-1} (-1)^{k+j}a_{j,k}
\left|
   \begin{smallmatrix}
    a_{21}      & \cdots & a_{2,k-1}   & a_{2,k+1}   & \cdots & a_{2,n-1}   \\
    \vdots      &        & \vdots      & \vdots      &        &   \vdots    \\
    a_{n-1,1}   & \cdots & a_{n-1,k-1} & a_{n-1,k+1} & \cdots & a_{n-1,n-1} \\
   \end{smallmatrix}
\right| \right) 
\end{multline}
But the parenthesis on the right-hand side is exactly expansion along
the $j$'th row of the following determinants
\begin{equation} \label{matrix3}
  \begin{cases}
   \left|
     \begin{smallmatrix}
      a_{11} & \cdots & a_{1,n-1} \\
     \vdots  &        & \vdots    \\
      a_{n-1,1} & \cdots & a_{n-1,n-1}
     \end{smallmatrix}
\right|,
& j=1 \\
  \left|
     \begin{smallmatrix}
       a_{21} & \cdots & a_{2,n-1} \\
      \vdots  &        & \vdots    \\
      a_{j,1} & \cdots & a_{j,n-1}\\
      a_{j,1} & \cdots & a_{j,n-1}\\
      \vdots  &        & \vdots    \\
      a_{n,1} & \cdots & a_{n,n-1}\\
     \end{smallmatrix}
\right|=0,   & 2\leq j \leq n-1 \\
-  \left|
     \begin{smallmatrix}
      a_{21} & \cdots & a_{2,n-1} \\
     \vdots  &        & \vdots    \\
      a_{n,1} & \cdots & a_{n,n-1}
     \end{smallmatrix}
\right|,   & j=n.
  \end{cases}
\end{equation} 
Combining (\ref{matrix2}) and (\ref{matrix3}) we obtain the right-hand
side of (\ref{matrix1}) and thus also the proof of (\ref{Hankela}).

To prove (\ref{Hankelb}) we prove the more general equation
\begin{multline} \label{matrix4}
 \left|
   \begin{smallmatrix}
    a_{21} & a_{22}    & \cdots & a_{2,n} \\
    a_{31} & a_{32}    & \cdots & a_{3,n} \\
    a_{41} & a_{42}    & \cdots & a_{4,n} \\
    \vdots & \vdots    &        & \vdots    \\
 a_{n+1,1} & a_{n+1,2} & \cdots & a_{n+1,n}
   \end{smallmatrix}
\right| \left|
   \begin{smallmatrix}
    a_{12} & a_{13}   & \cdots & a_{1,n}   \\
    a_{22} & a_{23}   & \cdots & a_{2,n}   \\
    \vdots & \vdots   &        & \vdots    \\
   a_{n-1,2} & a_{n-1,3}  & \cdots & a_{n-1,n}
   \end{smallmatrix}
\right| \\ = 
 \left|\begin{smallmatrix}
    a_{11} & a_{12}    & \cdots & a_{1,n} \\
    a_{21} & a_{22}    & \cdots & a_{2,n} \\
    \vdots & \vdots    &        & \vdots    \\
    a_{n-1,1} & a_{n-1,2}    & \cdots & a_{n-1,n} \\
 a_{n+1,1} & a_{n+1,2} & \cdots & a_{n+1,n}
   \end{smallmatrix}
\right|  \left|
   \begin{smallmatrix}
    a_{22} & a_{23}   & \cdots & a_{2,n}   \\
    a_{32} & a_{33}   & \cdots & a_{3,n}   \\
    \vdots & \vdots   &        & \vdots    \\
   a_{n,2} & a_{n,3}  & \cdots & a_{n,n}
   \end{smallmatrix}
\right| \\- 
 \left|\begin{smallmatrix}
    a_{11} & a_{12}    & \cdots & a_{1,n} \\
    a_{21} & a_{22}    & \cdots & a_{2,n} \\
    \vdots & \vdots    &        & \vdots    \\
    a_{n-1,1} & a_{n-1,2}    & \cdots & a_{n-1,n} \\
 a_{n,1} & a_{n,2} & \cdots & a_{n,n}
   \end{smallmatrix}
\right|  \left|
   \begin{smallmatrix}
    a_{22} & a_{23}   & \cdots & a_{2,n}   \\
    a_{32} & a_{33}   & \cdots & a_{3,n}   \\
    \vdots & \vdots   &        & \vdots    \\
   a_{n-1,2} & a_{n-1,3}  & \cdots & a_{n-1,n} \\
   a_{n+1,2} & a_{n+1,3}  & \cdots & a_{n+1,n} \\
   \end{smallmatrix}
\right|
\end{multline}
for $a_{ij}\in\mathbb C$, $i\in\{1,\ldots n+1\}$ and
$j\in\{1,\ldots,n\}$. We remark that Hankel-matrices are symmetric and for
these (\ref{matrix4}) reduces to (\ref{Hankelb}).
Observe that for $k\in \{2,\ldots,n\}$ we have
\begin{multline*}
0 = (-1)^k \left|
   \begin{smallmatrix}
a_{1,k}   &   a_{11}   & a_{12}     & \cdots & a_{1,n}   \\
a_{2,k}   &   a_{12}   & a_{22}     & \cdots & a_{2,n}   \\
\vdots    &   \vdots   & \vdots     &        & \vdots    \\
a_{n,k}   &  a_{n,2}   & a_{n,3}    & \cdots & a_{n,n} \\
a_{n+1,k} &  a_{n+1,2} & a_{n+1,3}  & \cdots & a_{n+1,n} \\
   \end{smallmatrix}
\right| \\ = (-1)^k\sum_{j=1}^{n+1}a_{j,k}(-1)^{j+1} 
 \left|
   \begin{smallmatrix}
  a_{11}      & a_{12}     & \cdots & a_{1,n}   \\
  \vdots      & \vdots     &        & \vdots    \\ 
  a_{j-1,1}   & a_{j-1,2}  & \cdots & a_{j-1,n}   \\
  a_{j+1,1}   & a_{j+1,2}  & \cdots & a_{j+1,n}   \\
  \vdots      & \vdots     &        & \vdots    \\
  a_{n+1,1} & a_{n+1,2}  & \cdots & a_{n+1,n} \\
   \end{smallmatrix}
\right|,
\end{multline*}
where the $j=1$ and $j=n+1$ are interpreted as remove the
1$^{\text{st}}$ and $(n+1)^{\text{th}}$ coloumn respectively. Thus also 
\begin{multline*}
0= \sum_{k=2}^n  \left|
   \begin{smallmatrix}
  a_{22}      & \cdots & a_{2,k-1} & a_{2,k+1} & \cdots & a_{2,n}   \\
  a_{32}      & \cdots & a_{3,k-1} & a_{3,k+1} & \cdots & a_{3,n}   \\
 \vdots      &        &\vdots     &  \vdots   &        & \vdots    \\ 
  a_{n-1,2}      & \cdots & a_{n-1,k-1} & a_{n-1,k+1} & \cdots & a_{n-1,n}   \\
   \end{smallmatrix}
\right| \\
\cdot
\left((-1)^k\sum_{j=1}^{n+1}a_{j,k}(-1)^{j+1} 
 \left|
   \begin{smallmatrix}
  a_{11}      & a_{12}     & \cdots & a_{1,n}   \\
  \vdots      & \vdots     &        & \vdots    \\ 
  a_{j-1,1}   & a_{j-1,2}  & \cdots & a_{j-1,n}   \\
  a_{j+1,1}   & a_{j+1,2}  & \cdots & a_{j+1,n}   \\
  \vdots      & \vdots     &        & \vdots    \\
  a_{n+1,1} & a_{n+1,2}  & \cdots & a_{n+1,n} \\
   \end{smallmatrix}
\right|\right)
\end{multline*}
Switching the indices we have
\begin{multline}\label{matrix5}
0=\sum_{j=1}^{n+1}  \left|
   \begin{smallmatrix}
  a_{11}      & a_{12}     & \cdots & a_{1,n}   \\
  \vdots      & \vdots     &        & \vdots    \\ 
  a_{j-1,1}   & a_{j-1,2}  & \cdots & a_{j-1,n}   \\
  a_{j+1,1}   & a_{j+1,2}  & \cdots & a_{j+1,n}   \\
  \vdots      & \vdots     &        & \vdots    \\
  a_{n+1,1} & a_{n+1,2}  & \cdots & a_{n+1,n} \\
   \end{smallmatrix}
\right| \\
\cdot \left(\sum_{k=2}^n (-1)^{k+j-1} a_{j,k} 
 \left|
   \begin{smallmatrix}
  a_{22}      & \cdots & a_{2,k-1} & a_{2,k+1} & \cdots & a_{2,n}   \\
  a_{32}      & \cdots & a_{3,k-1} & a_{3,k+1} & \cdots & a_{3,n}   \\
 \vdots      &        &\vdots     &  \vdots   &        & \vdots    \\ 
  a_{n-1,2}      & \cdots & a_{n-1,k-1} & a_{n-1,k+1} & \cdots & a_{n-1,n}   \\
   \end{smallmatrix}
\right|
\right) 
\end{multline}
The parenthesis of (\ref{matrix5}) is expansion along the
$j^{\text{th}}$ row of the following expression except for $j=n+1$ where
we expand along the $n^{\text{th}}$ row.
\begin{equation} \label{matrix6}
  \begin{cases}
    \left|
   \begin{smallmatrix}
    a_{12} & a_{13}   & \cdots & a_{1,n}   \\
    a_{22} & a_{23}   & \cdots & a_{2,n}   \\
    \vdots & \vdots   &        & \vdots    \\
   a_{n-1,2} & a_{n-1,3}  & \cdots & a_{n-1,n} 
   \end{smallmatrix}
\right|   , & j=1 \\
0, & j\in\{2,\ldots, n-1\} \\
 \left|
   \begin{smallmatrix}
    a_{22} & a_{23}   & \cdots & a_{2,n}   \\
    a_{32} & a_{33}   & \cdots & a_{3,n}   \\
    \vdots & \vdots   &        & \vdots    \\
   a_{n,2} & a_{n,3}  & \cdots & a_{n,n}
   \end{smallmatrix}
\right|,  & j=n \\
 -\left|
   \begin{smallmatrix}
    a_{22} & a_{23}   & \cdots & a_{2,n}   \\
    \vdots & \vdots   &        & \vdots    \\
   a_{n-1,2} & a_{n-1,3}  & \cdots & a_{n-1,n}   \\
   a_{n+1,2} & a_{n+1,3}  & \cdots & a_{n+1,n}   \\
   \end{smallmatrix}
\right| & j=n+1.
  \end{cases}
\end{equation}
Combining (\ref{matrix5}) and (\ref{matrix6}) we obtain (\ref{matrix4})
and this finishes the proof of (\ref{Hankelb}).
\end{proof}

\end{document}